\documentclass[a4paper,11pt]{article}
\usepackage{latexsym}
\usepackage{amssymb}
\usepackage{enumerate}
\usepackage{amsfonts}
\usepackage{amsmath}
\usepackage[dvips]{graphicx}
\usepackage[paper=a4paper,left=30mm,right=20mm,top=25mm,bottom=30mm]{geometry}
    \newenvironment{dedication}
        {\vspace{2ex}\begin{quotation}\begin{center}\begin{em}}
        {\par\end{em}\end{center}\end{quotation}}
\newcommand{\qed}{$\Box$}

\newenvironment{@abssec}[1]{%
    \if@twocolumn

      \section*{#1}%
    \else

      \vspace{.05in}\footnotesize
      \parindent .2in
 {\upshape\bfseries #1. }\ignorespaces
    \fi}

    {\if@twocolumn\else\par\vspace{.1in}\fi}

\newenvironment{keywords}{\begin{@abssec}{\keywordsname}}{\end{@abssec}}

\newenvironment{AMS}{\begin{@abssec}{\AMSname}}{\end{@abssec}}

\newcommand\keywordsname{Key words}
\newcommand\AMSname{AMS subject classifications}
\newcommand\AMname{AMS subject classification}
\newtheorem{theorem}{Theorem}
 \newtheorem{lemma}[theorem]{Lemma}
 
 \newtheorem{proposition}[theorem]{Proposition}
\newtheorem{remark}[theorem]{Remark}
 
\def\qed{\vbox{\hrule height0.6pt\hbox{%
  \vrule height1.3ex width0.6pt\hskip0.8ex
  \vrule width0.6pt}\hrule height0.6pt
 }}
\def\theequation{\arabic{section}.\arabic{equation}}
 \def\thetheorem{\arabic{section}.\arabic{theorem}}

\def\theequation{\arabic{section}.\arabic{equation}}
 \def\thetheorem{\arabic{section}.\arabic{theorem}}

\title{Some characterizations of parallel hyperplanes  in multi-layered heat conductors
\thanks{This research was partially supported by the Grants-in-Aid
for Scientific Research (B) ($\sharp$ 18H01126 and $\sharp$ 17H02847) of
Japan Society for the Promotion of Science.}}

\author{Shigeru Sakaguchi\thanks{Research Center for Pure and Applied Mathematics,
Graduate School of  Information Sciences, Tohoku
University, Sendai, 980-8579,  Japan.
({\tt sigersak@tohoku.ac.jp}).}}

\date{}
\begin{document}

\maketitle

\begin{dedication}
\begin{large}
Dedicated to Masaru Ikehata on the occasion of his 60th birthday
\end{large}
\end{dedication}

\begin{abstract}
We consider the Cauchy problem for the heat diffusion equation in the whole space consisting of three layers with different constant conductivities, where initially the upper and middle layers have temperature 0 and the lower layer has temperature 1. Under some appropriate conditions, it is shown that, if either the interface between the lower layer and the middle layer  is a stationary isothermic surface or there is a stationary isothermic surface in the middle layer near the lower layer, then the two interfaces must be parallel hyperplanes.  Similar propositions hold true,  either if a stationary isothermic surface is replaced by  a surface with the constant flow property or if the Cauchy problem is replaced by an appropriate initial-boundary value problem.

\vskip 2ex
\centerline{\bf R\'esum\'e}
\vskip 1ex

Nous consid\'erons le probl\`eme de Cauchy pour l'\'equation de diffusion de la chaleur dans tout l'espace compos\'e de trois couches avec diff\'erentes conductivit\'es constantes, o\`u initialement les couches sup\'erieure et moyenne ont la temp\'erature 0 et la couche inf\'erieure a la temp\'erature 1. Dans certaines conditions appropri\'ees, il est montr\'e que, si l'interface entre la couche inf\'erieure et la couche interm\'ediaire est une surface isotherme stationnaire ou s'il existe une surface isothermique stationnaire dans la couche interm\'ediaire pr\`es de la couche inf\'erieure, alors les deux interfaces doivent \^etre des hyperplans parall\`eles. Des propositions similaires sont vraies, soit si une surface isotherme stationnaire est remplac\'ee par une surface avec la propri\'et\'e d'\'ecoulement constant ou si le probl\`eme de Cauchy est remplac\'e par un probl\`eme de valeur de limite initiale appropri\'e.
  \end{abstract}



\pagestyle{plain}
\thispagestyle{plain}

\section{Introduction}
\label{introduction}
For $x \in \mathbb R^N$ with $N \ge 2$,  set $x = (x_1,\dots, x_{N-1}, x_N) = (y, x_N)$ for $y \in \mathbb R^{N-1}$. Let $f, h \in C^2(\mathbb R^{N-1})$ satisfy
$$
f(y)  < h(y) \ \mbox{ for every } y \in \mathbb R^{N-1}.
$$
Define two domains $D, \Omega$ in $\mathbb R^N$ by
\begin{equation}
\label{D and Omega}
D = \{ x \in \mathbb R^N\ :\ x_N > h(y) \},\quad \Omega =\{ x \in \mathbb R^N\ :\ x_N > f(y) \},
\end{equation}
respectively. Denote by $\sigma=\sigma(x)\ (x \in \mathbb R^N)$  the conductivity distribution of the whole medium given by
\begin{equation}
\label{conductivity constants}
\sigma =
\begin{cases}
\sigma_c \quad&\mbox{in } D, \\
\sigma_s \quad&\mbox{in } \Omega \setminus D, \\
\sigma_m \quad &\mbox{in } \mathbb R^N \setminus \Omega,
\end{cases}
\end{equation}
where $\sigma_c, \sigma_s, \sigma_m$ are positive constants with $\sigma_c  \not= \sigma_s$. This kind of three-phase electrical conductor has been dealt with in \cite{KLS2016} in the study of neutrally coated inclusions.

Let $u =u(x,t)$ be the unique bounded solution of either the Cauchy problem for the heat diffusion equation:
\begin{equation}
  u_t =\mbox{ div}(\sigma \nabla u)\quad\mbox{ in }\  \mathbb R^N\times (0,+\infty) \ \mbox{ and }\ u\ ={\mathcal X}_{\Omega^c}\ \mbox{ on } \mathbb R^N\times
\{0\},\label{heat Cauchy}
\end{equation}
where ${\mathcal X}_{\Omega^c}$ denotes the characteristic function of the set $\Omega^c=\mathbb R^N \setminus\Omega$, 
or the initial-boundary value problem for the heat diffusion equation:
\begin{eqnarray}
&&u_t =\mbox{ div}( \sigma \nabla u)\quad\mbox{ in }\ \Omega\times (0,+\infty), \label{heat equation initial-boundary}
\\
&&u=1  \ \quad\qquad\qquad\mbox{ on } \partial\Omega\times (0,+\infty), \label{heat Dirichlet}
\\ 
&&u=0  \  \quad\qquad\qquad \mbox{ on } \Omega\times \{0\}.\label{heat initial}
\end{eqnarray}
Let $g \in C^0(\mathbb R^{N-1})$ satisfy
$$
f(y)  < g(y) <  h(y) \ \mbox{ for every } y \in \mathbb R^{N-1}.
$$
Consider a domain $G$ in $\mathbb R^N$  defined by
\begin{equation}
\label{a testing domain G}
G = \{ x \in \mathbb R^N\ :\ x_N > g(y) \}.
\end{equation}
Suppose that 
\begin{equation}
\label{close to the interface below}
\mbox{ dist}(x, \partial\Omega) \le \mbox{ dist}(x, \overline{D}) \ \mbox{ for every } x \in \partial G.
\end{equation}
This assumption is  technical and corresponds to those in \cite[(5)]{Strieste2016},  \cite[(1.5)]{SBessatsu2017} and \cite[(1.6)]{CMSarXiv2018}, and it enables us to utilize the balance laws \cite[Theorem 2.1 and Corollary 2.2]{MSannals2002}.

Let us first state two theorems concerning stationary isothermic surfaces.
\begin{theorem}
\label{th:stationary isothermic1} Either let $N \le 8$ or let $\nabla f $ be bounded in $\mathbb R^{N-1}$ with $N\ge2$.  Suppose that $\partial\Omega$ is uniformly of class $C^6$ and the function $h-f$ has a minimum value in $\mathbb R^{N-1}$ and
moreover,  either  $h-f$ has a maximum value in $\mathbb R^{N-1}$ or  $h-f$ is unbounded in $\mathbb R^{N-1}$.
Let $u$ be the solution of problem \eqref{heat Cauchy}. If there exists a function $a : (0, +\infty) \to (0, +\infty) $ satisfying
\begin{equation}
\label{stationary isothermic surface}
u(x,t) = a(t)\ \mbox{ for every } (x,t) \in \partial\Omega \times (0, +\infty),
\end{equation}
then $\partial\Omega$ and $\partial D$ must be parallel hyperplanes.
\end{theorem}

\begin{theorem}
\label{th:stationary isothermic2} Either let $N \le 3$ or  let $\{ |f(y)-f(\hat{y})| :  |y-\hat{y}| \le 1 \}$ be bounded.  Suppose that  the function $h-f$ has a minimum value in $\mathbb R^{N-1}$ and
either  $h-f$ has a maximum value in $\mathbb R^{N-1}$ or  $h-f$ is unbounded in $\mathbb R^{N-1}$.  Let $u$ be the solution of problem \eqref{heat Cauchy} or problem \eqref{heat equation initial-boundary}-\eqref{heat initial}.  Under the assumption \eqref{close to the interface below},  if there exists a function $a : (0, +\infty) \to (0, +\infty) $ satisfying
\begin{equation}
\label{stationary isothermic surface near the interface below}
u(x,t) = a(t)\ \mbox{ for every } (x,t) \in \partial G \times (0, +\infty),
\end{equation}
then $\partial\Omega$ and $\partial D$ must be parallel hyperplanes.
\end{theorem}

In Theorems \ref{th:stationary isothermic1} and \ref{th:stationary isothermic2}, the conditions \eqref{stationary isothermic surface} and \eqref{stationary isothermic surface near the interface below} mean that each of $\partial\Omega$ and $\partial G$ is  a stationary isothermic surface. 
Thus each of  Theorems \ref{th:stationary isothermic1} and \ref{th:stationary isothermic2} characterizes parallel hyperplanes as the interfaces in such a way that there exists a stationary isothermic surface in the multi-layered heat conductors.  The assumptions on the function $h-f$ are technical, and in particular the existence of its maximum value or its minimum value enables us to utilize Hopf's boundary point lemma.

Next two theorems replace a stationary isothermic surface by  a surface with the constant flow property  which was dealt with in \cite{CMSarXiv2018}.

\begin{theorem}
\label{th:constant flow1} Either let $N \le 8$ or let $\nabla f $ be bounded in $\mathbb R^{N-1}$ with $N \ge 2$.  Suppose that $\partial\Omega$ is uniformly of class $C^6$ and the function $h-f$ has a minimum value in $\mathbb R^{N-1}$ and
moreover,  either  $h-f$ has a maximum value in $\mathbb R^{N-1}$ or  $h-f$ is unbounded in $\mathbb R^{N-1}$.
Let $u$ be the solution of problem \eqref{heat equation initial-boundary}-\eqref{heat initial}. If there exists a function $b : (0, +\infty) \to \mathbb R $ satisfying
\begin{equation}
\label{surface of the constant flow property}
\sigma_s \frac {\partial u}{\partial\nu}(x,t) = b(t)\ \mbox{ for every } (x,t) \in \partial\Omega \times (0, +\infty),
\end{equation}
then $\partial\Omega$ and $\partial D$ must be parallel hyperplanes, where $\nu$ denotes the outward unit normal vector to $\partial\Omega$.
\end{theorem}

\begin{theorem}
\label{th:constant flow2} Either let $N \le 3$ or let $\{ |f(y)-f(\hat{y})| :  |y-\hat{y}| \le 1 \}$ be bounded.  Suppose that  the function $h-f$ has a minimum value in $\mathbb R^{N-1}$ and
 either  $h-f$ has a maximum value in $\mathbb R^{N-1}$ or  $h-f$ is unbounded in $\mathbb R^{N-1}$,  and  moreover $g \in C^1(\mathbb R^{N-1})$.  Let $u$ be the solution of problem \eqref{heat Cauchy} or problem \eqref{heat equation initial-boundary}-\eqref{heat initial}.  Under the assumption \eqref{close to the interface below},  if there exists a function $b : (0, +\infty) \to \mathbb R $ satisfying
\begin{equation}
\label{surface of the constant flow property near the interface below}
\sigma_s \frac {\partial u}{\partial\nu}(x,t) = b(t)\ \mbox{ for every } (x,t) \in \partial G \times (0, +\infty),
\end{equation}
then $\partial\Omega$ and $\partial D$ must be parallel hyperplanes, where $\nu$ denotes the outward unit normal vector to $\partial G$.
\end{theorem}

 In Theorem \ref{th:constant flow1} the condition  \eqref{surface of the constant flow property},  together with the boundary condition \eqref{heat Dirichlet},  is overdetermined and it implies that the heat flow is parallel to the normal vector to $\partial \Omega$ and the amount of the flow is constant on $\partial\Omega$ for each time.  Such a condition was given by \cite{AG1989ARMA, GS2001pams}  for parabolic problems, which generalizes the overdetermined condition of Serrin \cite{Se1971ARMA} for elliptic problems.  Recently  such a boundary $\partial\Omega$ was called a surface with the constant flow property in the context of the heat flow in smooth Riemannian manifolds by  \cite{Sav2016MA}.  The condition \eqref{surface of the constant flow property near the interface below}, which  was introduced by  \cite{CMSarXiv2018},  is an overdetermination  different from Serrin-type, and we still  called it  the constant flow property in  \cite{CMSarXiv2018}.
 Similar characterizations of concentric balls in multi-phase heat conductors were obtained in the previous papers \cite{Strieste2016, SBessatsu2017, CMSarXiv2018}, and in the present paper we deal with  hyperplanes, which are not compact and need additional cares. The proofs of all the theorems consist of two steps. In the first step we show that $\partial\Omega$ must be a hyperplane, and the second step is devoted to proving that $\partial D$ is a hyperplane parallel to $\partial\Omega$.  We have two strategies in the first step; one applies to Theorems \ref{th:stationary isothermic2} and \ref{th:constant flow2} and the other  does to Theorems \ref{th:stationary isothermic1} and \ref{th:constant flow1}.  On the other hand,  the second step follows from one strategy common to all the theorems, which depends on a result concerning an elliptic overdetermined problem (see Theorem \ref{th:elliptic overdetermined} in section \ref{section5}).

\vskip 2ex
The following sections are organized as follows. In section \ref{section2}, we recall one lemma and three propositions from \cite{CMSarXiv2018,  Strieste2016},  where we need to modify the two propositions in order to deal with the case where $\partial\Omega$ is unbounded.  Indeed, we show that our case is reduced to the case where  $\partial\Omega$ is bounded and of class $C^2$ with the aid of the maximum principle and the Gaussian bounds for the fundamental solution of $u_t=\mbox{ div}(\sigma\nabla u)$ due to Aronson \cite[Theorem 1, p. 891]{Ar1967bams}(see also \cite[p. 328]{FaS1986arma}). Section \ref{section3} is devoted to the proofs of Theorems  \ref{th:stationary isothermic2} and \ref{th:constant flow2}; the balance laws (Proposition \ref{balance law}) and the asymptotic formula of the heat content of balls touching at a point on $\partial\Omega$ (Proposition \ref{modified asymptotic formula}) play a key role to show that $\partial\Omega$ must be a sort of Weingarten surface, and hence some results of \cite{SIndAM2013} implies that $\partial\Omega$ is a hyperplane.  Finally,  by using  Theorem 
\ref{th:elliptic overdetermined} given in section \ref{section5}, which concerns an elliptic overdetermined problem, we complete the proofs through the Laplace transform.
Section \ref{section4} is devoted to the proofs of Theorems \ref{th:constant flow1} and \ref{th:stationary isothermic1}. Under the assumption that $\partial\Omega$ is uniformly of class $C^6$,  the same arguments with the precise barriers as in the proofs of \cite[Theorems 1.4 and 1.5 in section 5]{CMSarXiv2018} work and we conclude that the mean curvature of $\partial\Omega$ must be constant even if $\partial\Omega$ is unbounded. Hence both the Bernstein theorem  and Moser's theorem for the minimal surface equation imply that $\partial\Omega$ is a hyperplane under appropriate assumptions. Finally,  Theorem \ref{th:elliptic overdetermined}  completes  the proofs through the Laplace transform. In section \ref{section5}, we give a proof of Theorem \ref{th:elliptic overdetermined}, where Hopf's boundary point lemma and the transmission condition on $\partial D$, together with three comparison principles and one  maximum principle for elliptic equations with discontinuous conductivities given in section \ref{section6},  play a key role.  Roughly, Theorem \ref{th:elliptic overdetermined} states that if $\partial\Omega$ is a hyperplane then $\partial D$ must be a hyperplane parallel to $\partial\Omega$. The last section \ref{section6} is devoted to the proofs of three comparison principles and one maximum principle for elliptic equations with discontinuous conductivities. 

\setcounter{equation}{0}
\setcounter{theorem}{0}

\section{Preliminaries}
\label{section2}
Let us introduce the distance function $\delta = \delta(x)$ of $x \in \mathbb R^N$  to $\partial\Omega$ by
\begin{equation}
\label{distance function to the boundary of the domain}
\delta(x) = \mbox{ dist}(x,\partial\Omega)\ \mbox{ for }\ x \in \mathbb R^N.
\end{equation}
We quote a lemma concerning the solutions of problem \eqref{heat Cauchy} and problem \eqref{heat equation initial-boundary}-\eqref{heat initial}
from \cite[Lemma 4.1]{CMSarXiv2018}, which simply comes from the maximum principle and the Gaussian bounds for the fundamental solution of $u_t=\mbox{ div}(\sigma\nabla u)$ due to Aronson \cite[Theorem 1, p. 891]{Ar1967bams}(see also \cite[p. 328]{FaS1986arma}).  Although \cite[Lemma 4.1]{CMSarXiv2018} concerns the case where $\Omega$ is bounded,  exactly the same proof is applicable  even if $\Omega$ is unbounded. 
For $\tau > 0$, we set
$$
\Omega_\tau = \{ x \in \Omega\ :\delta(x) \ge \tau \}\ \mbox{ and }\ \Omega_\tau^c
=\{ x \in \mathbb R^N\setminus\Omega\ : \delta(x) \ge \tau \}.
$$
\begin{lemma} 
\label{le:initial behavior and decay at infinity} 
Let $u$ be the solution of either problem \eqref{heat Cauchy} or problem \eqref{heat equation initial-boundary}-\eqref{heat initial} with a general conductivity
$\sigma=\sigma(x)\ (x\in \mathbb R^N)$ satisfying
$$
0 < \mu \le \sigma(x) \le M\ \mbox{ for every } x \in \mathbb R^N,
$$
where $\mu, M$ are positive constants. Then the following propositions hold true:
\begin{itemize}
\item[\rm (1)]  The solution $u$ satisfies
\begin{equation}
\label{between zero and one}
0 < u < 1\  \mbox{  in } \mathbb R^N \times (0,+\infty)\ \mbox{ or  in } \Omega \times (0,+\infty), \ \mbox{ respectively.}
\end{equation}
\item[\rm (2)] For every $\tau > 0$, there exist two positive constants $B$ and $b$ such that
$$
0 < u(x,t) < B e^{-\frac bt}\quad \mbox{ for every } (x,t) \in \Omega_\tau \times (0,+\infty)
$$
and, moreover, if $u$ is the solution of \eqref{heat Cauchy}, then
$$
0 < 1-u(x,t) < B e^{-\frac bt}\quad \mbox{ for every } (x,t) \in \Omega_\tau^c \times (0,+\infty).
$$
 \item[\rm(3)] The solution $u$ of \eqref{heat Cauchy} is such that 
 $$
 \lim\limits_{x \not\in \Omega, \delta(x) \to \infty} (1-u(x,t)) = 0  \quad \mbox{ for every } t \in (0,+\infty).
 $$
 \end{itemize}
\end{lemma}

In \cite[Theorems 1.3 and 1.2]{CMSarXiv2018}, a proposition (\cite[Proposition 2.2, pp. 171--172]{Strieste2016}) plays a key role, where the boundary of the domain is compact.
Here, we deal with the case where $\partial\Omega$  is unbounded, and therefore we need to modify the proposition. Denote by $B_r(x)$ an open ball in $\mathbb R^N$ with a radius $r >0$ and centered at a point $x \in \mathbb R^N$. The modified one is the following: 

\begin{proposition}
\label{modified asymptotic formula}
Let $\Omega$ be a possibly unbounded domain in $\mathbb R^N$, and let $x_0 \in \Omega$ and $z_0 \in \partial\Omega$.  Assume that $B_r(x_0) \subset \Omega, \ \overline{B_r(x_0)}\cap\partial\Omega =\{z_0\}$ and there exists  $\varepsilon > 0$ such that $\partial\Omega\cap B_\varepsilon(z_0)$ is of class $C^2$ and $\partial\Omega$ divides $B_\varepsilon(z_0)$ into two connected components.  Let $\sigma=\sigma(x)\ (x\in \mathbb R^N)$ be 
a general conductivity satisfying
$$
0 < \mu \le \sigma(x) \le M\ \mbox{ for every } x \in \mathbb R^N, \mbox{ and }\sigma(x) = \begin{cases} \sigma_s &\mbox{ if } x \in B_\varepsilon(z_0) \cap \Omega,\\  \sigma_m \ &\mbox{ if } x \in B_\varepsilon(z_0) \setminus \Omega,
\end{cases}
$$
where $\mu, M, \sigma_s, $ and $\sigma_m$ are positive constants. Let $u$ be the bounded solution of either problem \eqref{heat Cauchy} or problem \eqref{heat equation initial-boundary}-\eqref{heat initial} for this general conductivity $\sigma$.  Then we have:
\begin{equation}
\label{asymptotics and curvatures}
\lim_{t\to +0}t^{-\frac{N+1}4 }\!\!\!\!\!\int\limits_{B_r(x_0)}\!\!\! u(x,t)\, dx=
C(N, \sigma)\left\{\prod\limits_{j=1}^{N-1}\left(\frac 1r - \kappa_j(z_0)\right)\right\}^{-\frac 12}.
\end{equation}
Here, 
$\kappa_1(z_0),\dots,\kappa_{N-1}(z_0)$ denote the principal curvatures of $\partial\Omega$ at $z_0$ with 
respect to the inward normal direction to $\partial\Omega$  
and $C(N, \sigma)$ is a positive constant given by
$$
C(N,\sigma) = \left\{\begin{array}{rll}2\sigma_s^{\frac {N+1}4}c(N) \ &\mbox{ for problem \eqref{heat equation initial-boundary}-\eqref{heat initial} },
\\
\frac {2\sqrt{\sigma_m}}{\sqrt{\sigma_s}+\sqrt{\sigma_m}}\sigma_s^{\frac {N+1}4}c(N) &\mbox{ for problem \eqref{heat Cauchy} },
\end{array}\right.
$$
where $c(N)$ is a positive constant depending only on $N$. (Notice that if $\sigma_s=\sigma_m$ then $C(N, \sigma) = \sigma_s^{\frac {N+1}4}c(N)$
for problem \eqref{heat Cauchy}, that is, just half of the constant  for problem \eqref{heat equation initial-boundary}-\eqref{heat initial}.)
When $\kappa_j(z_0) = 1/r$ for some $j \in \{ 1, \cdots, N-1\}$, 
\eqref{asymptotics and curvatures} holds by setting the right-hand side to $+\infty$ (notice that  $\kappa_j(z_0) \le 1/r$ always holds for all $j$'s).
\end{proposition} 

\noindent
{\it Proof.\ }  It suffices to show that our case is reduced to the case where $\partial\Omega$ is bounded and of class $C^2$. Since  $\partial\Omega\cap B_\varepsilon(z_0)$ is of class $C^2$,  we can find a bounded domain $\Omega_*$  with $C^2$ boundary $\partial\Omega_*$ satisfying
$$
B_r(x_0)\cup \left(\Omega\cap \overline{B_{\frac 23\varepsilon}(z_0)}\right) \subset \Omega_*\subset \Omega,\  \overline{B_{\frac 23\varepsilon}(z_0)}\cap\partial\Omega\subset \partial\Omega_* 
\mbox{ and } \overline{B_r(x_0)}\cap\partial\Omega_* = \{z_0\}.
$$
Let us first consider problem  \eqref{heat equation initial-boundary}-\eqref{heat initial}.   Let $u_*=u_*(x,t)$ be the bounded solution of problem \eqref{heat equation initial-boundary}-\eqref{heat initial} where $\Omega$ and $\sigma$ are replaced with $\Omega_*$ and $\sigma_s$,  respectively. Then, it follows from \cite[Proposition 2.2, pp. 171--172]{Strieste2016} that the formula \eqref{asymptotics and curvatures} holds true for $u_*$.  We observe that the difference $v=u - u_*$ satisfies
\begin{eqnarray}
&&v_t =\sigma_s\Delta v\quad\mbox{ in }\ \left(\Omega\cap B_{\frac 23\varepsilon}(z_0)\right)\times (0,+\infty), \label{heat equation initial-boundary*}
\\
&&v=0  \ \quad\qquad\qquad\mbox{ on } \left(\partial\Omega\cap \overline{B_{\frac 23\varepsilon}(z_0)} \right)\times (0,+\infty), \label{heat Dirichlet*-1}
\\ 
&& |v| < 1\  \ \qquad\qquad\mbox{ on } \Omega_* \times (0,+\infty), \label{heat Dirichlet*-2}
\\
&&v=0  \  \quad\qquad\qquad \mbox{ on } \Omega_*\times \{0\}.\label{heat initial*}
\end{eqnarray}
Set
$$
\mathcal N = \left\{ x \in \mathbb R^N : \mbox{ dist}(x, \Omega_*\cap \partial B_{\frac 23\varepsilon}(z_0)) < \frac 1{100}\varepsilon \right\}.
$$
By comparing $v$ with the solutions of the Cauchy problem for the heat equation with conductivity $\sigma_s$ and initial data $\pm 2\mathcal X_{\mathcal N}$ for a short time, we see that
there exist two positive constants $B$ and $b$ such that
\begin{equation}
\label{exponential decay of the difference near the touching point}
|v(x,t)| \le Be^{-\frac bt}\ \mbox{ for every } (x,t) \in \overline{B_{\frac 12\varepsilon}(z_0) \cap \Omega} \times (0,\infty).
\end{equation}
By (2) of Lemma \ref{le:initial behavior and decay at infinity}, we may also have
\begin{equation}
\label{exponential decay inside the domains}
0 < u(x,t), \ u_*(x,t)  \le Be^{-\frac bt}\ \mbox{ for every } (x,t) \in \left(\overline{B_r(x_0)}\setminus B_{\frac 12\varepsilon}(z_0) \right)\times (0,\infty).
\end{equation}
Then,  it follows from \eqref{exponential decay of the difference near the touching point} and \eqref{exponential decay inside the domains}  that $u$ also satisfies  \eqref{asymptotics and curvatures}, since we already know that $u_*$ satisfies \eqref{asymptotics and curvatures}. Indeed,  observing that 
$$
t^{-\frac {N+1}4}\!\!\!\int\limits_{B_r(x_0)}\!\!\!v\, dx =t^{-\frac {N+1}4}\!\!\!\!\!\!\!\!\!\!\!\! \int\limits_{B_r(x_0)\setminus B_{\frac 12\varepsilon}(z_0)}\!\!\!\!\!\!\!\!\!\!\!\!v\, dx  + t^{-\frac {N+1}4}\!\!\!\!\!\!\!\!\!\!\!\!\int\limits_{B_r(x_0)\cap B_{\frac 12\varepsilon}(z_0)}\!\!\!\!\!\!\!\!\!\!\!\!v\, dx 
$$
and letting $t \to \infty$ yield the conclusion.

It remains to consider problem \eqref{heat Cauchy}.  Let us define the conductivity $\sigma_*=\sigma_*(x)\ (x \in \mathbb R^N)$ by
\begin{equation}
\label{conductivity constants two-phase}
\sigma_* =
\begin{cases}
\sigma_s \quad&\mbox{in } \Omega_*, \\
\sigma_m \quad &\mbox{in } \mathbb R^N \setminus \Omega_*.
\end{cases}
\end{equation}
Let $u_*=u_*(x,t)$ be the bounded solution of problem \eqref{heat Cauchy} where $\Omega$ and $\sigma$ are replaced with $\Omega_*$ and $\sigma_*$,  respectively.
Then, it follows from \cite[Proposition 2.2, pp. 171--172]{Strieste2016} that the formula \eqref{asymptotics and curvatures} holds true for $u_*$.  
We observe that the difference $v=u - u_*$ satisfies
\begin{eqnarray}
&&v_t =\mbox{ div}(\sigma_*\nabla v)\quad\mbox{in }\ B_{\frac 23\varepsilon}(z_0)\times (0,+\infty), \label{heat equation initial-boundary*C}
\\ 
&& |v| < 1\  \ \qquad\qquad\mbox{ in }\ \mathbb R^N \times (0,+\infty), \label{heat bounds-2C}
\\
&&v=0  \  \quad\qquad\qquad \mbox{ on } \left(\Omega_*\cup  \overline{B_{\frac 23\varepsilon}(z_0)}\right)\times \{0\}.\label{heat initial*C}
\end{eqnarray}
Then, by the same comparison arguments with the aid of the Gaussian bounds due to Aronson \cite[Theorem 1, p. 891]{Ar1967bams}(see also \cite[p. 328]{FaS1986arma}), we see that  there exist two positive constants $B$ and $b$ satisfying \eqref{exponential decay inside the domains} and 
\begin{equation}
\label{exponential decay of the difference near the touching point C}
|v(x,t)| \le Be^{-\frac bt}\ \mbox{ for every } (x,t) \in \overline{B_{\frac 12\varepsilon}(z_0)} \times (0,\infty),
\end{equation}
and hence $u$ also satisfies  \eqref{asymptotics and curvatures}. \qed

 Since a proposition \cite[Proposition E]{CMSarXiv2018}, where the boundary of the domain is compact,  also plays a key role in \cite{CMSarXiv2018},  
 we need to modify the proposition in order to deal with the case where $\partial\Omega$ is unbounded.
 
  \begin{proposition}
\label{modified formula on the boundary value}
Let $\Omega$ be a possibly unbounded domain in $\mathbb R^N$, and let $z_0 \in \partial\Omega$.  Assume that  there exists  $\varepsilon > 0$ such that $\partial\Omega\cap B_\varepsilon(z_0)$ is of class $C^2$ and $\partial\Omega$ divides $B_\varepsilon(z_0)$ into two connected components.  Let $\sigma=\sigma(x)\ (x\in \mathbb R^N)$ be 
a general conductivity satisfying
$$
0 < \mu \le \sigma(x) \le M\ \mbox{ for every } x \in \mathbb R^N, \mbox{ and }\sigma(x) = \begin{cases} \sigma_s &\mbox{ if } x \in B_\varepsilon(z_0) \cap \Omega,\\  \sigma_m \ &\mbox{ if } x \in B_\varepsilon(z_0) \setminus \Omega,
\end{cases}
$$
where $\mu, M, \sigma_s, $ and $\sigma_m$ are positive constants. Let $u$ be the bounded solution of  problem \eqref{heat Cauchy} for this general conductivity $\sigma$.  Then, as $t \to +0$, $u$ converges to the number $\frac {\sqrt{\sigma_m}}{\sqrt{\sigma_s}+\sqrt{\sigma_m}}$ uniformly on $\partial\Omega\cap \overline{B_{\frac 12\varepsilon}(z_0)}$.
\end{proposition} 

\noindent
{\it Proof.\ } It suffices to show that our case is reduced to the case where $\partial\Omega$ is bounded and of class $C^2$. 
As in the proof of Proposition \ref{modified asymptotic formula} for  problem \eqref{heat Cauchy},  let $u_*=u_*(x,t)$ be the bounded solution of problem \eqref{heat Cauchy} where $\Omega$ and $\sigma$ are replaced with $\Omega_*$ and $\sigma_*$,  respectively. Then $u_*$ satisfies the conclusion because of \cite[Proposition E]{CMSarXiv2018}. Therefore,  since $v = u-u_*$ satisfies \eqref{exponential decay of the difference near the touching point C},  $u$ also satisfies the conclusion. \qed

We quote another ingredient called a balance law adjusted to our use  from \cite[Lemma 4.2]{CMSarXiv2018} and \cite[Theorem 2.1]{MSannals2002}. 
For convenience, we give a proof with the aid of \cite[Theorem 2.1]{MSannals2002}.


 \begin{proposition}[\cite{CMSarXiv2018, MSannals2002}]
\label{balance law}
Let $W$ be a domain in $\mathbb R^N$ with $N \ge 2$, and let $u = u(x,t)$ satisfy
$$
u_t = \sigma_s\Delta u\ \mbox{ in } W \times (0, +\infty).
$$
Consider two points $p, q \in W$ and two unit vectors $\xi, \eta \in \mathbb R^N$.  Set 
$$
r_*=\min\{\mbox{\rm dist}(p, \partial W), \mbox{\rm dist}(q, \partial W)\}.
$$
Then the following three propositions hold true:
\begin{itemize}
\item[\rm (1)]\ $u(p,t) = u(q,t)$ for every $t > 0$ if and only if 
$$
\int\limits_{B_r(p)}\!\!u(x,t)\, dx = \int\limits_{B_r(q)}\!\!u(x,t)\, dx \ \mbox{ for every }(t,r ) \in (0,+\infty)\times(0,r_*).
$$
\item[\rm (2)]\ $\xi\cdot\nabla u(p,t) = \eta\cdot\nabla u(q,t)$ for every $t > 0$ if and only if 
$$
\xi\cdot\!\!\!\!\int\limits_{B_r(p)}\!\!u(x,t) (x-p)\, dx = \eta\cdot\!\!\!\!\int\limits_{B_r(q)}\!\!u(x,t)(x-q)\, dx \ \mbox{ for every }(t,r ) \in (0,+\infty)\times(0,r_*).
$$
\item[\rm (3)]\ $\nabla u(p,t) = 0$ for every $t > 0$ if and only if 
$$
\int\limits_{B_r(p)}\!\! u(x,t)(x-p)\, dx = 0 \ \mbox{ for every }(t,r ) \in (0,+\infty)\times(0, \mbox{\rm dist}(p, \partial W)).
$$
\end{itemize}
\end{proposition}

\noindent
{\it Proof.\ } (3) is just \cite[Corollary 2.2]{MSannals2002}.   (1) follows from \cite[Theorem 2.1]{MSannals2002}. Indeed, consider the function 
$$
v_1(x,t) = u(x+p,t) - u(x+q,t)\ \mbox{ for } (x,t) \in B_{r_*}(0) \times (0,+\infty).
$$
Then $v_1$ satisfies the heat equation with conductivity $\sigma_s$ and $v_1(0,t) = 0$ for every $t > 0$. Thus  \cite[Theorem 2.1]{MSannals2002} gives the conclusion.

(2) is proved in  \cite[Lemma 4.2]{CMSarXiv2018} with the aid of  \cite[Theorem 2.1]{MSannals2002}.  For (2), by choosing an orthogonal matrix $A$ satisfying $A\xi = \eta$, we consider the function
$$
v_2(x,t) = u(x+p,t) - u(Ax+q,t)\ \mbox{ for } (x,t) \in B_{r_*}(0) \times (0,+\infty).
$$
Then the function $\xi\cdot \nabla v_2(x,t)$ satisfies the heat equation with conductivity $\sigma_s$ and for every $t > 0$
$$
\xi\cdot \nabla v_2(0,t) = \xi\cdot \nabla u(p, t) -\eta\cdot \nabla u(q,t) =0.
$$ 
Thus, it follows from \cite[Theorem 2.1]{MSannals2002} that 
$$
\xi\cdot\!\!\!\!\int\limits_{B_r(0)}\nabla v_2(x,t)\, dx = 0 \ \mbox{ for every } (t,r ) \in (0,+\infty)\times(0,r_*),
$$
and hence, by the divergence theorem and again integrating in $r$, we infer that 
$$
\xi\cdot\!\!\!\!\int\limits_{B_r(0)}v_2(x,t)x\, dx = 0\ \mbox{ for every } (t,r ) \in (0,+\infty)\times(0,r_*),
$$
which gives (2). \qed

\setcounter{equation}{0}
\setcounter{theorem}{0}

\section{Proofs of Theorems \ref{th:stationary isothermic2} and \ref{th:constant flow2}:  the 1st strategy}
\label{section3}
Under each of the assumptions of Theorems \ref{th:stationary isothermic2} and \ref{th:constant flow2}, we follow the proofs of \cite[Theorems 1.1 and 1.3]{Strieste2016} and \cite[Theorem 1.2]{CMSarXiv2018}, respectively,  in order to prove that $\partial\Omega$ is parallel to $\partial G$ and the quantity $\displaystyle \prod_{j=1}^{N-1}(1/R-\kappa_j(z))$ is constant for $z \in \partial\Omega$, where $R$ is the distance between $\partial\Omega$ and  $\partial G$,  $\kappa_1(z), \dots, \kappa_{N-1}(z)$ denote the principal curvatures of $\partial\Omega$ at a point $z \in \partial\Omega$ with respect to the inward normal direction to $\partial\Omega$, and $\displaystyle\max_{1\le j \le N-1}\kappa_j < 1/R$  for every $z \in \partial\Omega$.  Once this is proved, we immediately infer that  $\partial\Omega$ must be a hyperplane. Indeed,  if $N=2$ then $\partial\Omega$ must be a straight line,  if $N=3$, by \cite[Theorem 4, p. 281]{SIndAM2013},  $\partial\Omega$ must be a hyperplane, and  if $\{ |f(y)-f(\hat{y})| :  |y-\hat{y}| \le 1 \}$ is bounded with $N\ge 2$, by \cite[Theorem 3 and Remark 3, p. 273]{SIndAM2013},  the same conclusion holds true.  In  the proof of \cite[Theorem 4, p. 281]{SIndAM2013}, the strong comparison principle for the viscosity solutions of the minimal surface equation plays a key role.  Note that  \cite{OSNonlAnal2019} gives a simple proof of the strong comparison principle for the prescribed mean curvature equation including the minimal surface equation. 

We need to modify \cite[Lemma 4.3]{CMSarXiv2018} in order to deal with the case where $\partial\Omega$ is unbounded  and $\partial G$ is of class $C^1$ under the assumption \eqref{surface of the constant flow property near the interface below}. 
\begin{lemma} 
\label{le: constant weingarten curvature}
Let $u$ be the solution of either problem \eqref{heat Cauchy} or problem \eqref{heat equation initial-boundary}--\eqref{heat initial}. 
Under each of the assumptions \eqref{stationary isothermic surface near the interface below} and \eqref{surface of the constant flow property near the interface below} of {\rm Theorems \ref{th:stationary isothermic2} and \ref{th:constant flow2}},  the following assertions hold:
\begin{enumerate}[\rm (1)]
\item there exists a number $R > 0$ such that 
$$
\delta(x)= R\ \mbox{ for every } x \in \partial G,
$$
where $\delta(x)$ is the distance function given by \eqref{distance function to the boundary of the domain};
\item $\partial\Omega$ and $\partial G$ are real analytic hypersurfaces;
\item  the mapping $\partial\Omega \ni z \mapsto x(z)  \equiv z-R\,\nu(z) \in \partial G$ is a diffeomorphism where $\nu(z)$ denotes the outward unit normal vector to $\partial\Omega$ at $z \in \partial\Omega$; in particular $\partial\Omega$ and $\partial G$  are parallel hypersurfaces at distance $R$;
\item  the principal curvatures of $\partial\Omega$ satisfy
\begin{equation*}
\label{bounds of curvatures}
 \max_{1\le j\le N-1}\kappa_j(z) < \frac 1{R}\ \mbox{ for every } z \in \partial\Omega;
\end{equation*}
\item there exists a number $c > 0$ satisfying 
\begin{equation}
\label{monge-ampere}
\prod_{j=1}^{N-1}\left(\frac 1R-\kappa_j(z)\right) = c\ \mbox{ for every } z \in \partial\Omega.
\end{equation}
\end{enumerate}
\end{lemma}

Before proving this lemma, we prepare a purely geometric lemma for the proof of  Theorem \ref{th:constant flow2}.
\begin{lemma} 
\label{le: purely geometric}
Suppose that $g \in C^1(\mathbb R^{N-1})$ in the definition \eqref{a testing domain G} of $G$.  Set
 $$
 R = \inf\{ \delta(x)\ :\ x \in \partial G \} ( \ge 0),
 $$
where $\delta(x)$ is the distance function given by \eqref{distance function to the boundary of the domain}.  Then, for every $\varepsilon > 0$, there exists a point
$p \in\partial G$ such that
\begin{eqnarray}
&&\delta(p) < R + \varepsilon; \label{near the infmum}
\\ 
&&\mbox{ there exists a point } z \in \partial\Omega \mbox{ with } B_{\delta(p)}(p) \cap \partial\Omega = \{ z\}; \label{unique projection point}
\\
&& (z-p)\cdot\nu(p) \not=0\ \mbox{ and } \max_{1\le j \le N-1}\kappa_j(z) < \frac 1{\delta(p)},\label{good point for the asymptotic formula}
\end{eqnarray}
where $\nu(p)$ denotes the outward unit normal vector to $\partial G$ at $p \in \partial G$.
\end{lemma}

\noindent
{\it Proof.\ } Let $\varepsilon > 0$. Set
$$
G_\varepsilon =\left\{ x \in \mathbb R^N\ :\ x_N > g(y)+\frac \varepsilon2 \right\}. 
$$
Since $\inf\{ \delta(x)\ :\ x \in \partial G_\varepsilon \} \le R+\frac \varepsilon2$,  there exists a point $q \in \partial G_\varepsilon$ with $ \delta(q) < R+\varepsilon$. Then there exists $z \in \partial\Omega$ with $\delta(q) = |q-z|$.  By the intermediate value theorem there exists a point $p \in 
\partial G \cap \overline{qz}$ such that
$$
|p-z| < |q-z| < R+\varepsilon,
$$
where $ \overline{qz} $ denotes the line segment connecting $q$ and $z$. Therefore we infer that
\begin{equation}
\label{finding good point}
B_{\delta(p)}(p) \cap \partial\Omega = \{ z\}\ \mbox{ and }\ \max_{1\le j \le N-1} \kappa_j(z) \le \frac 1{\delta(q)} < \frac 1{\delta(p)}.
\end{equation}
Hence, by  the inverse mapping theorem and \eqref{finding good point},  there exists an infinite solid cylinder $U$, whose axis is the line containing  $\overline{qz}$, such that
$$
\delta \in C^2(\overline{U \cap (\Omega\setminus G)})\ \mbox{ and } \nabla \delta(p) = \frac {p-z}{|p-z|}.
$$
If $\nabla\delta(p)\cdot\nu(p) \not=0$, then the conclusion follows from \eqref{finding good point}. Thus, let us consider the case where $\nabla\delta(x)\cdot\nu(x) = 0$ for all $x \in U \cap (\Omega\setminus G) \cap \partial G$. Let $x = x(s)\ (s \in \mathbb R)$ the curve determined by the Cauchy problem:
\begin{equation}
\label{Cauchy problem for the curve}
\frac d{ds} x(s) = -\nabla\delta(x(s))\ \mbox{ and } x(0) = p.
\end{equation}
Then, as long as $x(s)$ exists,  $x(s) \in \partial G$ and moreover, since $\nabla\delta(x) = \frac {p-z}{|p-z|}$ for every $x \in  \overline{pz}$,  we have from the uniqueness of the solution of the Cauchy problem \eqref{Cauchy problem for the curve}
$$
x(s) = p - s\frac {p-z}{|p-z|}.
$$
 These contradict the fact that  $\delta(x(s)) \ge R$ and $\delta(x(s)) = -s +\delta(p)$. Thus there exists a point $x \in  U \cap (\Omega\setminus G) \cap \partial G$ with $\nabla\delta(x)\cdot\nu(x) \not= 0.$ This point $x \in \partial G$ replaces $p$. \qed

\vskip 2ex
\noindent
{\bf Proof of Lemma \ref{le: constant weingarten curvature}.\ } First, it follows from the assumption \eqref{close to the interface below} that 
$$
B_r(x) \subset \Omega\setminus\overline{D} \ \mbox{ for every } x \in \partial G\ \mbox{ with  } 0 < r \le \delta(x).
$$
Therefore, since $\sigma = \sigma_s$ in $\Omega\setminus\overline{D}$, we can use  Lemma \ref{balance law}. 

Let us first deal with Theorem \ref{th:stationary isothermic2}. Then, with the aid of Lemma \ref{balance law}, Lemma \ref{le:initial behavior and decay at infinity} and Proposition \ref{modified asymptotic formula},  under the assumption \eqref{stationary isothermic surface near the interface below} of  Theorem \ref{th:stationary isothermic2} the same proof as in \cite[Lemma 2.4, pp. 176--179]{Strieste2016} is applicable in showing all the assertions (1)--(5) of this lemma even if $\partial\Omega$ is not compact. Roughly,  suppose that 
$\delta(p)< \delta(q)$ for some points $p, q \in \partial G$.  Then, \eqref{stationary isothermic surface near the interface below} gives (1) of Proposition \ref{balance law}. In particular, we choose $r = \delta(p)$. On the other hand, combining  (2) of Lemma \ref{le:initial behavior and decay at infinity}  and Proposition \ref{modified asymptotic formula} yields a contradiction to (1) of Proposition \ref{balance law} with $r = \delta(p)$. Thus assertion (1) holds under  the assumption \eqref{stationary isothermic surface near the interface below}.  Once we have (1) under the assumption \eqref{stationary isothermic surface near the interface below} of  Theorem \ref{th:stationary isothermic2}, the others (2)--(5) follow easily.  In particular, the analyticity of $\partial G$ follows from the analyticity of the solution $u = u(x,t)$ in $x$, if one shows that for every $x\in \partial G$ there exists a time $t > 0$ satisfying $\nabla u(x,t) \not=0$ with the aid of  \eqref{stationary isothermic surface near the interface below},  (3) of Lemma \ref{balance law}, (2) of Lemma \ref{le:initial behavior and decay at infinity} and Proposition \ref{modified asymptotic formula}.  $\partial \Omega$ is also real analytic by (3).

Let us proceed to Theorem \ref{th:constant flow2}.  Since  \cite[Lemma 4.3]{CMSarXiv2018} concerns the case where $\partial\Omega$  is compact and $\partial G$ is of class $C^2$, we need to modify its proof  in order to deal with the case where $\partial\Omega$ is not compact  and $\partial G$ is of class $C^1$. Let us consider assertion (1) under  the assumption \eqref{surface of the constant flow property near the interface below} of Theorem \ref{th:constant flow2}.
Let $\varepsilon > 0$. Then it follows from Lemma \ref{le: purely geometric} that there exists a point $p \in \partial G$ satisfying \eqref{near the infmum}--\eqref{good point for the asymptotic formula}. Hence it follows from Proposition \ref{modified asymptotic formula}  and  (2) of Lemma \ref{le:initial behavior and decay at infinity} that 
\begin{equation}
\label{asymptotics and curvatures with constant flow}
\lim_{t\to +0}t^{-\frac{N+1}4 }\!\nu(p)\cdot\!\!\!\!\!\!\!\int\limits_{B_{\delta(p)}(p)}\!\!\!\!\!\!\! u(x,t)(x-p)\, dx=
C(N, \!\sigma)\nu(p)\!\cdot\!(z-p)\!\left\{\prod\limits_{j=1}^{N-1}\!\left(\frac 1{\delta(p)}\! - \!\kappa_j(z)\right)\right\}^{-\frac 12}\!\!\!\!\!\! \not= 0.
\end{equation}
Suppose that there exists a point $q \in \partial G$ with $\delta(p) < \delta(q)$.  Then, \eqref{surface of the constant flow property near the interface below} gives (2) of Proposition \ref{balance law}. In particular, we choose $r = \delta(p),\ \xi = \nu(p)$ and $ \eta =\nu(q)$ to infer that
\begin{equation}
\label{balance for constant heat flow}
t^{-\frac{N+1}4 }\nu(p)\cdot\!\!\!\!\!\!\!\!\int\limits_{B_{\delta(p)}(p)}\!\!\!\!\!\!\!u(x,t) (x-p)\, dx =t^{-\frac{N+1}4 } \nu(q)\cdot\!\!\!\!\!\!\!\!\int\limits_{B_{\delta(p)}(q)}\!\!\!\!\!\!\!u(x,t)(x-q)\, dx \ \mbox{ for every } t > 0.
\end{equation}
 On the other hand,  it follows from (2) of Lemma \ref{le:initial behavior and decay at infinity} that the right-hand side of \eqref{balance for constant heat flow} tends to $0$ as $t \to +0$, which contradicts \eqref{asymptotics and curvatures with constant flow}.  Therefore, we conclude that $\delta(q) \le \delta(p)$ for every $q \in \partial G$.
Moreover, \eqref{near the infmum} yields that $\delta(q) = R$  for every $q \in \partial G$ and $R > 0$. Thus assertion (1) holds also under  the assumption \eqref{surface of the constant flow property near the interface below}. 

Once we have (1) under the assumption \eqref{surface of the constant flow property near the interface below} of Theorem \ref{th:constant flow2}, we infer that for every $x \in \partial G$ there exists a unique $z =z(x) \in \partial \Omega$ satisfying
\begin{equation}
\label{uniqueness of touching points}
\overline{B_R(x)}\cap\partial\Omega =\{z(x)\},
\end{equation}
since $\partial G$ is of class $C^1$.  As in \cite[Lemma 2.4, pp. 176--179]{Strieste2016}, we introduce the set $\gamma \subset \partial\Omega$ by
$$
\gamma = \{ z \in \partial\Omega : \overline{B_R(x(z))}\cap\partial\Omega =\{z\}\ \mbox{ for } x(z) =z-R\nu(z) \in \partial G\ \mbox{ and } \max_{1\le j\le N-1}\kappa_j(z) < 1/R \}.
$$
Then Lemma \ref{le: purely geometric} implies that $\gamma \not=\emptyset$, and assertion (1) yields that 
$$
B_R(z) \cap G = \emptyset\ \mbox{ and } \nu(x(z)) = \nu(z)  \mbox{ for every } z \in \gamma.
$$
Thus, we infer that  the formula \eqref{asymptotics and curvatures with constant flow} holds if we set $p=x(z) \in \partial G$  with $z \in \gamma$ and $\nu(p)\cdot(z-p) = R = \delta(p)$, that is, for every $z \in \gamma$
\begin{equation}
\label{asymptotics and curvatures with constant flow 2nd}
\lim_{t\to +0}t^{-\frac{N+1}4 }\!\nu(x(z))\cdot\!\!\!\!\!\!\!\!\!\int\limits_{B_{R}(x(z))}\!\!\!\!\!\!\! u(x,t)(x-x(z))\, dx=
C(N, \!\sigma)R\!\left\{\prod\limits_{j=1}^{N-1}\!\left(\frac 1{R}\! - \!\kappa_j(z)\right)\right\}^{-\frac 12} > 0.
\end{equation}
Hence, combining (2) of Proposition \ref{balance law} with this formula \eqref{asymptotics and curvatures with constant flow 2nd} yields that there exists a number $c > 0$ satisfying 
\begin{equation}
\label{monge-ampere pre}
\prod_{j=1}^{N-1}\left(\frac 1R-\kappa_j(z)\right) = c\ \mbox{ for every } z \in \gamma.
\end{equation}
Then, since $\partial\Omega$ is of class $C^2$,  combining \eqref{uniqueness of touching points} with \eqref{monge-ampere pre} yields that $\gamma$ is closed in $\partial\Omega$. On the other hand, the inverse mapping theorem implies that $\gamma$ is also open in $\partial\Omega$ and the mapping $\gamma \ni z \mapsto x(z) \in \partial G$ is a local diffeomorphism.  Therefore $\gamma = \partial\Omega$, since $\partial\Omega$ is connected. Thus the others (3)--(5) follow immediately.  Finally, the analyticity of $\partial\Omega$ follows from (5) and hence $\partial G$ is also real analytic by (3).  The proof of Lemma \ref{le: constant weingarten curvature} is completed. \qed 

\noindent
{\bf Completion of the proofs  of  Theorems \ref{th:stationary isothermic2} and \ref{th:constant flow2} : }
As mentioned in the beginning of this section,  Lemma \ref{le: constant weingarten curvature} implies that $\partial\Omega$ must be a hyperplane under each of the assumptions of Theorems \ref{th:stationary isothermic2} and \ref{th:constant flow2}.  Then, by Lemma \ref{le: constant weingarten curvature}, $\partial G$ must be a hyperplane parallel to $\partial\Omega$.  Let us prove Theorems \ref{th:stationary isothermic2} and \ref{th:constant flow2} by using Theorem \ref{th:elliptic overdetermined} given in section \ref{section5}. 

Let $u$ be the solution of problem \eqref{heat Cauchy}.
We introduce the function $w=w(x)\ (x \in \overline{\mathbb R^N})$ by
\begin{equation}
\label{Laplace transform with lambda = 1}
w(x) = \int_0^\infty\!\! e^{-t}u(x,t)\, dt.
\end{equation}
Then $w$ satisfies 
\begin{eqnarray}
&& -\mbox{ div}( \sigma \nabla w) + w = 0\ \mbox{ in } \Omega, \label{equation inside Omega}
\\
&& - \sigma_m \Delta (1-w) + (1-w)= 0\ \mbox{ in } \mathbb R^N \setminus \overline{\Omega}, \label{equation outside Omega}
\\
&& w|_- = w|_+ \mbox{ and } \sigma_s \frac {\partial w}{\partial\nu}|_- = \sigma_m \frac {\partial w}{\partial\nu}|_+ \ \mbox{ on } \partial\Omega,  \label{transmission on the boundary of Omega}
\\
&& 0 < w < 1\ \mbox{ in } \mathbb R^N, \label{all value in between 0 and 1}
\\
&& \lim_{x \not\in \Omega, \delta(x) \to \infty} (1-w(x)) = 0, \label{limit at infinity outside Omega}
\end{eqnarray}
where $+$ denotes the limit from outside and $-$ that from inside of $\Omega$ and \eqref{limit at infinity outside Omega} comes from (3) of Lemma \ref{le:initial behavior and decay at infinity} and Lebesgue's dominated convergence theorem. Then \eqref{equation inside Omega} and \eqref{all value in between 0 and 1} give \eqref{elliptic equation 1} and \eqref{between 0 and 1} in section \ref{section5}, respectively.  Thus it suffices to show \eqref{overdetermination usual}.
 Let $\Theta \in \mathbb R^N$ be an arbitrary vector parallel to the hyperplanes $\partial\Omega$ and $\partial G$. Consider the function
$$
v^*(x,t) = u(x,t) - u(x+\Theta,t)\  \mbox{ for } (x,t) \in \left(\mathbb R^N \setminus G\right) \times (0,+\infty).
$$
Then $v^* = v^*(x,t)$ satisfies
\begin{eqnarray*}
&&v^*_t =\mbox{ div}( \sigma \nabla v^*)\quad\mbox{ in }\ \left(\mathbb R^N \setminus \overline{G}\right) \times (0,+\infty), 
\\
&&\mbox{Either } v^*=0  \ \mbox{ on } \partial G\times (0,+\infty)\  \mbox{ or }\  \frac {\partial v^*}{\partial \nu} = 0 \ \mbox{ on } \partial G\times (0,+\infty), 
\\ 
&&v^*=0  \  \quad\qquad\qquad \mbox{ on }  \left(\mathbb R^N \setminus \overline{G}\right)\times \{0\}.
\end{eqnarray*}
Hence it follows from the maximum principle that $v^* \equiv 0$, that is,  in $\left(\mathbb R^N \setminus G\right) \times (0,+\infty)$, the solution $u$ depends only on $\delta(x)$ and $t$ since $\Theta \in \mathbb R^N$ is an arbitrary vector parallel to the hyperplane $\partial\Omega$. Therefore $w$ depends only on $\delta(x)$ in $\mathbb R^N \setminus G$ and hence \eqref{overdetermination usual} holds true.  \eqref{all value in between 0 and 1} gives the fact that $0 < \alpha < 1$ in \eqref{overdetermination usual}, and \eqref{equation outside Omega}, \eqref{transmission on the boundary of Omega} and \eqref{limit at infinity outside Omega} yield that $\beta > 0$.  Indeed, by solving \eqref{equation outside Omega},  we get
$$
1-w(x) = c_0 \exp\left(-\frac{\delta(x)}{\sqrt{\sigma_m}}\right) \ \mbox{ for every } x \in \mathbb R^N \setminus \Omega,
$$
for some positive number $0 < c_0 < 1.$ This together with \eqref{transmission on the boundary of Omega} yields that $\beta > 0$. Therefore Theorem \ref{th:elliptic overdetermined} implies the conclusion of Theorems \ref{th:stationary isothermic2} and \ref{th:constant flow2} for problem \eqref{heat Cauchy}.

It remains to take care of the solution $u$ of problem \eqref{heat equation initial-boundary}-\eqref{heat initial}. We introduce the function $w=w(x)\ (x \in \overline{\Omega})$ by \eqref{Laplace transform with lambda = 1}.  Then $w$ satisfies 
\begin{eqnarray}
&& -\mbox{ div}( \sigma \nabla w) + w = 0\ \mbox{ in } \Omega, \label{equation inside Omega initial-boundary}
\\
&& 0 < w < 1\ \mbox{ in } \Omega, \label{all value in between 0 and 1 inside Omega}
\\
&& w = 1\ \mbox{ on } \partial\Omega. \label{boundary value equals 1}
\end{eqnarray}
Hence \eqref{equation inside Omega initial-boundary} and \eqref{all value in between 0 and 1 inside Omega} give \eqref{elliptic equation 1} and \eqref{between 0 and 1} in section \ref{section5}, respectively.  Thus it suffices to show \eqref{overdetermination usual}.
 Let $\Theta \in \mathbb R^N$ be an arbitrary vector parallel to the hyperplanes $\partial\Omega$ and $\partial G$. Consider the function
$$
v^*(x,t) = u(x,t) - u(x+\Theta,t)\  \mbox{ for } (x,t) \in \left(\overline{\Omega}\setminus G\right) \times (0,+\infty).
$$
Then $v^* = v^*(x,t)$ satisfies
\begin{eqnarray*}
&&v^*_t =\mbox{ div}( \sigma \nabla v^*)\quad\mbox{ in }\ \left(\Omega \setminus \overline{G}\right) \times (0,+\infty), 
\\
&&\mbox{Either } v^*=0  \ \mbox{ on } \partial G\times (0,+\infty)\  \mbox{ or }\  \frac {\partial v^*}{\partial \nu} = 0 \ \mbox{ on } \partial G\times (0,+\infty), 
\\
&& v^*=0  \ \mbox{ on } \left[\partial\Omega\times (0,+\infty)\right]\cup  \left[ \left(\Omega \setminus \overline{G}\right)\times \{0\}\right].
\end{eqnarray*}
Hence it follows from the maximum principle that $v^* \equiv 0$, that is,  in $\left(\Omega \setminus G\right) \times (0,+\infty)$, the solution $u$ depends only on $\delta(x)$ and $t$ since $\Theta \in \mathbb R^N$ is an arbitrary vector parallel to the hyperplane $\partial\Omega$. Therefore $w$ depends only on $\delta(x)$ in $\Omega \setminus G$ and hence \eqref{overdetermination usual} holds true. \eqref{boundary value equals 1} gives that $\alpha = 1$,  and  it follows from \eqref{all value in between 0 and 1 inside Omega}, \eqref{boundary value equals 1} and Hopf's boundary point lemma that $\beta > 0$.  Therefore Theorem \ref{th:elliptic overdetermined} implies the conclusion of Theorems \ref{th:stationary isothermic2} and \ref{th:constant flow2} for problem \eqref{heat equation initial-boundary}-\eqref{heat initial}. \qed

\setcounter{equation}{0}
\setcounter{theorem}{0}

\section{Proofs of Theorems \ref{th:constant flow1} and \ref{th:stationary isothermic1}:  the 2nd strategy}
\label{section4}
Under the assumptions of Theorems \ref{th:constant flow1} and \ref{th:stationary isothermic1}, we follow the proofs of \cite[Theorems 1.4 and 1.5 in section 5]{CMSarXiv2018} in order to prove that the mean curvature of $\partial\Omega$ is constant. Once this is proved, we immediately infer that  $\partial\Omega$ must be a hyperplane. Indeed, since $\partial \Omega$ is an entire graph over $\mathbb R^{N-1}$, the constant mean curvature must be zero and if $N=2$ then $\partial\Omega$ must be a straight line, if $3 \le N \le 8$, by the Bernstein theorem for the minimal surface equation (see \cite[Theorem 17.8, p. 208]{G1984}),  $\partial\Omega$ must be a hyperplane, and if $\nabla f$ is bounded in $\mathbb R^{N-1}$ with $N \ge 3$, by Moser's theorem \cite[Corollary, p. 591]{M1961} (see also \cite[Theorem 17.5, p. 205]{G1984}),  the same conclusion holds true. 

Since $\partial\Omega$ is uniformly of class $C^6$, there exists two positive numbers $r$ and $K$ such that, for every point $p \in \partial\Omega$, there exist an orthogonal coordinate system $z$ and a function $\varphi \in C^6(\mathbb R^{N-1})$ such that the $z_N$ coordinate axis lies in the inward normal direction to $\partial\Omega$ at $p$, the origin is located at $p$,\   $C^6$ norm of $\varphi$ in $\mathbb R^{N-1}$  is less than $K$,\  $\varphi(0)=0,\ \nabla \varphi(0) = 0$ and the set $B_r(p) \cap\Omega$  is written as in the $z$ coordinate system
$$
 \{ z \in B_r(0) : z_N > \varphi(z_1, \dots, z_{N-1}) \}.
$$
Since $\partial\Omega$ is uniformly of class $C^6$ as explained above,  by choosing a number $\delta_0 > 0$ sufficiently small and setting  
\begin{equation}\label{inner tubular neighborhood of Omega}
\mathcal N_- = \{ x \in \Omega\ :\ 0< \delta(x) < \delta_0 \}\ \mbox{ and }\ \mathcal N_+ = \{ x \in \mathbb R^N \setminus \overline{\Omega}\ :\ 0< \delta(x) < \delta_0 \},
\end{equation}
where $\delta(x)$ is the distance function given by \eqref{distance function to the boundary of the domain}, 
we see that
\begin{eqnarray}
&& \delta \in C^6(\overline{\mathcal N_\pm}), \ \sup\left\{ \left|\frac {\partial^\alpha \delta}{\partial x^\alpha}(x)\right| : x \in \overline{\mathcal N_\pm}, |\alpha| \le 6 \right\} < +\infty,\ \sigma = \begin{cases} \sigma_s &\mbox{ in } \mathcal N_-,\\  \sigma_m &\mbox{ in } \mathcal N_+,
\end{cases}\quad
 \label{c6 regularity}
\\
&&\mbox{ for every } x \in \overline{\mathcal N_\pm} \mbox{ there exists a unique }z = z(x) \in\partial\Omega \mbox{ with } \delta(x) = |x-z|, \qquad\label{the nearest point z from x}
\\
&& z(x) = x -\delta(x)\nabla\delta(x)\ \mbox{ for all } x \in \overline{\mathcal N_\pm}, \qquad\label{ the point z and distance from x} 
\\
&& \max_{1\le j \le N-1}|\kappa_j(z)| < \frac 1{2\delta_0}\ \mbox{ for every } z \in \partial\Omega, \qquad\label{upper bound of the curvatures on the boundary}
\end{eqnarray}
where  $\kappa_1(z), \dots, \kappa_{N-1}(z)$ denote the principal curvatures of $\partial\Omega$ at a point $z \in \partial\Omega$ with respect to the inward normal direction $-\nu(z)=\nabla\delta(z)$ to $\partial\Omega$ for $\delta \in C^6(\overline{\mathcal N_-})$.

As in the proofs of \cite[Theorems 1.4 and 1.5 in section 5]{CMSarXiv2018}, we introduce the function $w = w(x,\lambda)$ by
$$
w(x,\lambda) = \lambda\int_0^\infty\!\! e^{-\lambda t}u(x,t)\, dt \
\begin{cases}
\mbox{ for } (x,\lambda) \in \overline{\Omega}\times(0,+\infty) &\mbox{ in  problem \eqref{heat equation initial-boundary}-\eqref{heat initial}},
\\
\mbox{ for } (x,\lambda) \in \mathbb R^N\times(0,+\infty) &\mbox{ in  problem \eqref{heat Cauchy}}.
\end{cases}
$$
Although the difference between \cite[Theorems 1.4 and 1.5]{CMSarXiv2018} and Theorems  \ref{th:constant flow1} and \ref{th:stationary isothermic1} is such that the neighborhoods of $\partial\Omega$ is bounded in \cite[Theorems 1.4 and 1.5]{CMSarXiv2018} and they are unbounded in Theorems  \ref{th:constant flow1} and \ref{th:stationary isothermic1}, we have all the ingredients corresponding to those in  \cite[Theorems 1.4 and 1.5]{CMSarXiv2018}; the maximum principle (Proposition \ref{prop:maximum principle on unbounded domains}) enables us to use the comparison arguments on each of unbounded neighborhoods $\mathcal N_\pm$; (2) of Lemma \ref{le:initial behavior and decay at infinity} yields that $w(x,\lambda)$ and $1-w(x,\lambda)$ decay exponentially as $\lambda \to \infty$ on $\partial \mathcal N_-\setminus\partial\Omega$ and $\partial\mathcal N_+\setminus \partial\Omega$,  respectively;  Proposition \ref{modified formula on the boundary value} works for problem \eqref{heat Cauchy} even if $\partial\Omega$ is unbounded;  the situation \eqref{c6 regularity}--\eqref{upper bound of the curvatures on the boundary} coming from the fact that $\partial\Omega$ is uniformly of class $C^6$ enables us to construct the same precise barriers for $w$;  and moreover, by introducing an increasing sequence of bounded subdomains in each of $\mathcal N_\pm$ together with an increasing sequence of bounded harmonic functions on each of the subdomains,  we can construct a harmonic function $\psi = \psi(x)$,  as the limit of the sequence,  on each of $\mathcal N_\pm$  satisfying 
$$
\psi = 0 \ \mbox{ on } \partial\Omega,\  \psi = 2\ \mbox{ on } \partial\mathcal N_\pm \setminus \partial\Omega  \mbox{ and }  0 < \psi < 2\ \mbox{ in } \mathcal N_\pm,
$$
even if $\mathcal N_\pm$ is unbounded. This harmonic function $\psi$ was needed in constructing the precise barriers  in the proofs of \cite[Theorems 1.4 and 1.5]{CMSarXiv2018}.  Therefore, the same arguments  as in the proofs of \cite[Theorems 1.4 and 1.5 in section 5]{CMSarXiv2018} work and we conclude that the mean curvature of $\partial\Omega$ must be constant.
Thus, as mentioned in the beginning of this section, $\partial\Omega$ must be a hyperplane.  Hence, as in the proofs  of  Theorems \ref{th:stationary isothermic2} and \ref{th:constant flow2} in section \ref{section3}, we may infer that $w(x) = w(x,1)$ satisfies \eqref{elliptic equation 1}--\eqref{overdetermination usual} with
$0 < \alpha \le 1$ and $\beta > 0$. Therefore Theorem \ref{th:elliptic overdetermined} implies the conclusion of Theorems  \ref{th:constant flow1} and \ref{th:stationary isothermic1}. \qed


\setcounter{equation}{0}
\setcounter{theorem}{0}

\section{An elliptic overdetermined problem}
\label{section5}
In this section, we assume that $\partial\Omega$ is a hyperplane, that is,  $f$ is an affine function in \eqref{D and Omega}.
Moreover,  let us assume that there exists a function $w= w(x)\ (x \in \overline{\Omega})$ which satisfies  the following:
\begin{eqnarray}
&&-\mbox{ div}( \sigma \nabla w) + w = 0\ \mbox{ in } \Omega, \label{elliptic equation 1}
\\
&&0 < w < 1 \mbox{ in } \Omega, \label{between 0 and 1}
\\
&&w = \alpha \ \mbox{ and } \  \sigma_s\frac{\partial w}{\partial\nu}= \beta \ \mbox{ on } \partial\Omega,\label{overdetermination usual}
\end{eqnarray}
where  $\nu$ denotes the outward unit normal vector to $\partial\Omega$, $\sigma$ is given by  \eqref{conductivity constants} and
$\alpha,\ \beta$ are constants with $0 < \alpha \le 1,\  \beta > 0$, respectively. Define two functions $w_{\pm}$ by
$$
w_+(x) = w(x)\ \mbox{ for } x \in \overline{\Omega} \setminus D\ \mbox{ and }\ w_-(x) = w(x)\ \mbox{ for } x \in \overline{D}.
$$
Then the transmission condition for $w$ on $\partial D$ is written as
\begin{equation}
\label{transmission condition on the interface above}
w_+ = w_-\ \mbox{ and }\ \sigma_s\frac {\partial w_+}{\partial\nu} = \sigma_c\frac {\partial w_-}{\partial\nu}\ \mbox{ on } \partial D,
\end{equation}
where $\nu$ denotes the outward unit normal vector to $\partial D$.


\begin{theorem}
\label{th:elliptic overdetermined} Suppose that  the function $h-f$ has a minimum value in $\mathbb R^{N-1}$ and
 either  $h-f$ has a maximum value in $\mathbb R^{N-1}$ or  $h-f$ is unbounded in $\mathbb R^{N-1}$. Then
$\partial D$ must be a hyperplane parallel to $\partial\Omega$.
\end{theorem}

\begin{remark}
We basically follow the arguments in {\rm \cite{SBessatsu2017}} to prove this theorem.  However, the difference is such that {\rm \cite{SBessatsu2017}} concerns concentric balls and {\rm Theorem \ref{th:elliptic overdetermined}} does parallel hyperplanes; the former is compact and the latter is not compact. As mentioned in section \ref{introduction},  Hopf's boundary point lemma and the transmission condition \eqref{transmission condition on the interface above} on $\partial D$, together with three comparison principles and one  maximum principle for elliptic equations with discontinuous conductivities given in section \ref{section6},  play a key role. 
\end{remark}

\noindent
{\bf Proof of Theorem \ref{th:elliptic overdetermined}.\ }  Since $\partial\Omega$ is a hyperplane, by a translation and a rotation we may assume that in the new coordinate system $z$
$$
\Omega =\{ z \in \mathbb R^N\ :\ z_N > 0 \}.
$$
Then, with the aid of the uniqueness of the solutions of the Cauchy problem for elliptic equations, we see that $w_+$ must be a function of  one variable $\rho=z_N$ and $w_+ = w_+(\rho)$ satisfies
\begin{equation}
\label{ODE Laplace transformed}
-\sigma_s w_+^{\prime\prime}(\rho) + w_+(\rho) = 0 \ \mbox{ in } \Omega \setminus \overline{D}, \ w_+(0) = \alpha \mbox{ and } \sigma_s w_+^\prime(0) = - \beta.
\end{equation}
Moreover we extend $w_+$ as a unique solution of the above Cauchy problem in \eqref{ODE Laplace transformed}  for all $\rho = z_N$  with $z \in \mathbb R^N$ and
we have for some constants $c_1, c_2$ 
\begin{equation}
\label{ODE explicit solutions for sigma s}
w_+(\rho) = c_1 \exp\left(-\frac \rho{\sqrt{\sigma_s}}\right) +  c_2 \exp\left(\frac \rho{\sqrt{\sigma_s}}\right)   \mbox{ for all }\ \rho \in \mathbb R.
\end{equation}
Then it follows from \eqref{ODE Laplace transformed} that
$$
c_1+c_2 =\alpha \in (0, 1],\  \sqrt{\sigma_s}(c_1-c_2) = \beta > 0\ \mbox{ and hence }\ c_1 > 0.
$$
In view of the assumption, we may deal with the following two cases in the original coordinate system $x$:
 $$
 \mbox{ (I)  $h-f$ is unbounded in $\mathbb R^{N-1}$}; \quad  \mbox{ (II)  $h-f$ has a maximum value in $\mathbb R^{N-1}$}.
 $$
Let us consider case (I) first. \eqref{between 0 and 1} yields that $c_2 = 0$ and hence $0 < c_1 \le 1$ by \eqref{ODE Laplace transformed}. Thus
\begin{equation}
\label{decaying solution at infinity}
w_+(\rho) = c_1 \exp\left(-\frac \rho{\sqrt{\sigma_s}}\right) \ \mbox{ with } 0 < c_1 \le 1\  \mbox{ for all }\ \rho \in \mathbb R.
\end{equation}
Then we notice that 
\begin{equation}
\label{monotone decreasing}
w_+^\prime(\rho) < 0\ \mbox{ for all }\ \rho \in \mathbb R\ \mbox{ and }\ \lim_{\rho \to +\infty} w_+(\rho) = 0.
\end{equation}
Since the function $h-f$ has a minimum value in $\mathbb R^{N-1}$  and $f$ is an affine function in the original coordinate system $x$, there exists a point $z^* \in \partial D$ in the new coordinate system $z$ satisfying
$$
z^*_N = \min\limits_{z \in \partial D} z_N > 0\ \mbox{ and }\ \{ z_N \in \mathbb R : z \in \partial D \} = [z^*_N, \infty).
$$
Let $v_* = v_*(\rho)\ (\rho \ge z^*_N)$ be the unique solution of the Cauchy problem:
$$
-\sigma_c v_*^{\prime\prime}(\rho) +v_*(\rho) = 0\ \mbox{ for } \rho \in \mathbb R,\ v_*(z^*_N) = w_+(z^*_N)\ \mbox{ and }\ \sigma_c v_*^\prime(z^*_N) = \sigma_s w_+^\prime(z^*_N).
$$
 Hence we have for some constants $c^*_3, c^*_4$ 
\begin{equation}
\label{ODE explicit solutions for sigma c}
v_*(\rho) = c^*_3 \exp\left(-\frac \rho{\sqrt{\sigma_c}}\right) +  c^*_4 \exp\left(\frac \rho{\sqrt{\sigma_c}}\right)  \ \mbox{ for } \rho \in \mathbb R.
\end{equation}
Distinguish two cases:
$$
\mbox{ (I-1) } \ \sigma_c > \sigma_s; \quad \mbox{ (I-2) } \ \sigma_c < \sigma_s.
$$
In case (I-1) we have from \eqref{monotone decreasing} that
\begin{equation}
\label{transmission negative derivative}
\sigma_s w_+^\prime(z^*_N)  = \sigma_c  v_*^\prime(z^*_N) < 0.
\end{equation}
Hence,  with \eqref{monotone decreasing} in hand, by applying  (2)-(ii) of Proposition \ref{prop:ODE comparison with different conductivities} to $w_1= w_+$ and $w_2=v_*$,  we have
\begin{equation}
\label{comparison ode and c4}
 w_+(\rho) <  v_*(\rho)\ \mbox{ for every  }  \rho > z^*_N, \mbox{ and hence } c^*_4 \ge 0.
\end{equation}
We also have
\begin{equation}
\label{boundary condition on the boundary of D}
w \not\equiv v_* \mbox{ and } w \le v_* \ \mbox{ on } \partial D.
\end{equation}
Therefore, since $-\sigma_c \Delta w +w = -\sigma_c\Delta v_* + v_* = 0\ \mbox{ and } 0 < w < 1\mbox{ in } D,  \ c^*_4 \ge 0$ and $\min\{v_*,1\}$ is a bounded supersolution in $D$,  it follows from the comparison principle (Proposition \ref{prop:maximum principle on unbounded domains}) that
\begin{equation}
\label{inequality wanted in D from above}
v_* > w\ \mbox{ in } D.
\end{equation}
Here we applied Proposition \ref{prop:maximum principle on unbounded domains} to the function $\min\{v_*,1\}-w$ in $D$.
Thus, with the aid of Hopf's boundary point lemma at $z^* \in \partial D$,  this contradicts the fact that 
$$
v_* = w\ \mbox{ and }\ \frac {\partial v_*}{\partial \nu} = \frac {\partial w}{\partial \nu}\Big|_- \left(= \frac {\partial w_-}{\partial \nu}\right) \ \mbox{ at }\  z^*,  
$$
where $\nu$ denotes the outward unit normal vector to $\partial D$ and $-$ denotes the limit from inside of $D$. Here we used \eqref{transmission condition on the interface above}.

In case (I-2), we also have \eqref{transmission negative derivative} from \eqref{monotone decreasing} and the same argument as in case (I-1), together with Proposition \ref{prop:ODE comparison with different conductivities},  yields that \eqref{comparison ode and c4} is replaced with 
\begin{equation}
\label{comparison ode opposite and c4}
  v_*(\rho) <   w_+(\rho)\ \mbox{ for every  }  \rho > z^*_N, \mbox{ and hence } c^*_4 \le 0,
\end{equation}
and then the comparison principle (Proposition \ref{prop:maximum principle on unbounded domains}) gives
\begin{equation}
\label{inequality wanted in D from below}
v_* < w\ \mbox{ in } D,
\end{equation}
since $\max\{ v_*, 0 \}$ is a bounded subsolution in $D$.
Thus we get a contradiction with the aid of Hopf's boundary point lemma at $z^* \in \partial D$.  Therefore, case (I) does not occur.

Let us proceed to case (II).  Since the function $h-f$ has a maximum value in $\mathbb R^{N-1}$  and $f$ is an affine function in the original coordinate system $x$, there exists a point $z^\sharp \in \partial D$ in the new coordinate system $z$ satisfying
\begin{equation}
\label{the boundary of D in z coordinates}
z^\sharp_N = \max\limits_{z \in \partial D} z_N > 0\ \mbox{ and }\ \{ z_N \in \mathbb R : z \in \partial D \} = [z^*_N, z^\sharp_N].
\end{equation}
If $z^\sharp_N = z^*_N$, then $\partial D$ must be a hyperplane parallel to $\partial\Omega$ and hence the conclusion of Theorem \ref{th:elliptic overdetermined} holds true. Therefore we distinguish three cases:
$$
\mbox{ (IIa) } \ c_2 = 0\ \mbox{ and } z^\sharp_N > z^*_N ; \quad\mbox{ (IIb) } \ c_2 < 0\ \mbox{ and } z^\sharp_N > z^*_N ; \quad \mbox{ (IIc) } \  c_2 > 0 \ \mbox{ and } z^\sharp_N > z^*_N.
$$
In case (IIa)  \eqref{ODE explicit solutions for sigma s} yields \eqref{monotone decreasing}.
Then the same arguments as in case (I) work and we get a contradiction, that is, case (IIa) does not occur.

In case (IIb) we notice that \eqref{monotone decreasing} is replaced with
\begin{equation}
\label{monotone decreasing to -infinity}
w_+^\prime(\rho) < 0\ \mbox{ for all }\ \rho \in \mathbb R\ \mbox{ and }\ \lim_{\rho \to +\infty} w_+(\rho) = -\infty.
\end{equation}
Distinguish two cases:
$$
\mbox{ (IIb-1) } \ \sigma_c > \sigma_s; \quad \mbox{ (IIb-2) } \ \sigma_c < \sigma_s.
$$
With \eqref{monotone decreasing to -infinity} in hand, in case (IIb-2)  by the same arguments as in case (I-2) we notice that $c^*_4< 0$ and hence we obtain \eqref{inequality wanted in D from below} which gives a contradiction with the aid of Hopf's boundary point lemma at $z^* \in \partial D$. In case (IIb-1), if $c^*_4 \ge 0$, then the same  arguments as in case (I-1) also work and one can get a contradiction. Thus it suffices to take care of case (IIb-1) with $c^*_4 < 0$.

Let us consider case (IIb-1) with $c^*_4 < 0$. For every $r\ge z^*_N$, we introduce the solution $v_r = v_r(\rho)\ (\rho \in \mathbb R)$ of the Cauchy problem:
$$
-\sigma_c v_r^{\prime\prime}(\rho) +v_r(\rho) = 0\ \mbox{ for } \rho \in \mathbb R,\ v_r(r) = w_+(r)\ \mbox{ and }\ \sigma_c v_r^\prime(r) = \sigma_s w_+^\prime(r).
$$
 Hence we have for some constants $c_3(r), c_4(r)$ 
\begin{equation}
\label{ODE explicit solutions for sigma c with r}
v_r(\rho) = c_3(r) \exp\left(-\frac \rho{\sqrt{\sigma_c}}\right) +  c_4(r) \exp\left(\frac \rho{\sqrt{\sigma_c}}\right)  \ \mbox{ for } \rho \in \mathbb R.
\end{equation}
In particular, we have
\begin{equation}
\label{the function c_4 on r}
c_4(r) = \frac {\sqrt{\sigma_c} w_+(r)+\sigma_s w_+^\prime(r)}{2\sqrt{\sigma_c}}\exp\left(-\frac r{\sqrt{\sigma_c}}\right).
\end{equation}
Note that $c_3(z^*_N) = c^*_3,\ c_4(z^*_N) = c^*_4$ and $v_{z^*_N} = v_*$,  where $c^*_3,\ c^*_4$ and $v_*$ are given in \eqref{ODE explicit solutions for sigma c}. Set
\begin{equation}
\label{name sharp}
c^\sharp_3 = c_3(z^\sharp_N),\   c^\sharp_4 = c_4(z^\sharp_N) \mbox{ and }\ v_\sharp = v_{z^\sharp_N}.
\end{equation}
Distinguish two cases:
$$
\mbox{ (IIb-1-1) }\ c^\sharp_4 \le 0; \quad\mbox{ (IIb-1-2) }\ c^\sharp_4 > 0.
$$
In case (IIb-1-1),  with \eqref{monotone decreasing to -infinity} in hand,  the same arguments as in (I) also work and \eqref{comparison ode and c4} is replaced with 
\begin{equation}
\label{comparison ode and c4 sharp}
v_\sharp(\rho) < w_+(\rho) \ \mbox{ for every  }  \rho < z^\sharp_N.
\end{equation}
Then we also have
\begin{equation}
\label{boundary condition on the boundary of D at sharp}
w \not\equiv v_\sharp \mbox{ and } w \ge v_\sharp \ \mbox{ on } \partial D,
\end{equation}
and the comparison principle (Proposition \ref{prop:maximum principle on unbounded domains}) gives
\begin{equation}
\label{inequality wanted in D from below sharp}
v_\sharp < w\ \mbox{ in } D,
\end{equation}
since $\max\{ v_\sharp, 0 \}$ is a bounded subsolution in $D$.
Thus we get a contradiction with the aid of Hopf's boundary point lemma at $z^\sharp \in \partial D$.  Therefore, case  (IIb-1-1) does not occur.

In case (IIb-1-2),  in view of \eqref{monotone decreasing to -infinity} and \eqref{the function c_4 on r}, we observe that there exists $R > 0$ satisfying
$$
z^*_N < z^\sharp_N < R,\ c^*_4 = c_4(z^*_N) < 0, \ c^\sharp_4= c_4(z^\sharp_4) > 0\ \mbox{ and } c_4(R) < 0.
$$
By \eqref{the function c_4 on r}, $c_4(r)$ is continuous in $r$. Therefore, it follows from the intermediate value theorem that there exist two numbers $r_1$ and $r_2$ satisfying
$$
z^*_N < r_1 < z^\sharp_N < r_2 < R\ \mbox{ and }\ c_4(r_1) = c_4(r_2) = 0,
$$
and hence in particular both the functions $v_{r_j}\ (j=1, 2)$ are bounded in $[0,\infty)$.
Introduce  two functions $w_j = w_j(\rho)\ (j=1, 2)$ for $\rho \ge 0$ by
$$
w_j(\rho) = 
 \begin{cases}
w_+(\rho)  &\mbox{if }\  0 \le \rho \le r_j, \\
v_{r_j}(\rho) &\mbox{if }\ \rho >  r_j.
\end{cases} 
$$
Then we can apply Proposition \ref{prop:ODE comparison with discontinuous conductivities} to these $w_j = w_j(\rho)\ (j=1, 2)$ and obtain that
$r_1=r_2$, which is a contradiction. Therefore, case  (IIb-1-2) does not occur.

In case (IIc) it follows that there exists a unique $r_0 > 0$ satisfying
$$
w_+^\prime(\rho) < 0\ \mbox{ if }\ \rho < r_0,\  w_+^\prime(\rho) > 0\ \mbox{ if }\ \rho > r_0\ \mbox{ and } \lim_{\rho \to +\infty}w_+(\rho) = +\infty.
$$
Distinguish three cases:
$$
\mbox{ (IIc-1) }\ 0 < r_0 \le z^*_N; \quad\mbox{ (IIc-2) }\ z^*_N < r_0 < z^\sharp_N;  \quad\mbox{ (IIc-3) }\ z^\sharp_N \le r_0.
$$
Let us first consider case (IIc-1).  Distinguish two cases:
$$
\mbox{ (IIc-1-1) } \ \sigma_c > \sigma_s; \quad \mbox{ (IIc-1-2) } \ \sigma_c < \sigma_s.
$$
In case (IIc-1-1),  we employ $v_\sharp$.  It follows from (1) of Proposition \ref{prop:ODE comparison with different conductivities} that
\begin{equation}
\label{inequality on the boundary of D needed}
w_+(\rho)  < v_\sharp(\rho) \ \mbox{if } r_0 \le \rho < z_N^\sharp.
\end{equation}
Moreover,  by integrating the ordinary differential equations which $w_+$ and $v_\sharp$ satisfy, we have
$$
-\sigma_c v_\sharp^\prime(r_0)= -\left(\sigma_c v_\sharp^\prime(r_0) - \sigma_s w_+^\prime(r_0)\right) = \int_{r_0}^{z_n^\sharp}(v_\sharp(\rho)-w_+(\rho)) d\rho > 0.
$$
Hence we notice that 
$$
v_\sharp^\prime(r_0) < 0\ \mbox{ and  }\ v_\sharp^\prime(z_N^\sharp) = \frac {\sigma_s}{\sigma_c} w_+^\prime(z_N^\sharp) > 0. 
$$
This implies that $v_\sharp$ must have a critical point and hence $c_4^\sharp > 0$.   We  also have from \eqref{inequality on the boundary of D needed} that
\begin{equation}
\label{boundary condition on the boundary of D at sharp above}
w \not\equiv v_\sharp \mbox{ and } w \le v_\sharp \ \mbox{ on } \partial D.
\end{equation}
Thus the comparison principle (Proposition \ref{prop:maximum principle on unbounded domains}) gives
\begin{equation}
\label{inequality wanted in D from above sharp}
v_\sharp > w\ \mbox{ in } D,
\end{equation}
since $\min\{ v_\sharp, 1 \}$ is a bounded supersolution in $D$ because of the fact that $c_4^\sharp > 0$.
Thus we get a contradiction with the aid of Hopf's boundary point lemma at $z^\sharp \in \partial D$.  Therefore, case  (IIc-1-1) does not occur.

In case  (IIc-1-2),  we employ $v_*$ instead of $v_\sharp$.  It follows from (1) of Proposition \ref{prop:ODE comparison with different conductivities} that
\begin{equation}
\label{inequality on the boundary of D needed from above}
w_+(\rho)  < v_*(\rho) \ \mbox{if } \rho > z_N^*, \mbox{ and hence } c_4^* > 0.
\end{equation}
Here positivity of  $c_4^*$ comes from that of $c_2$.  Thus the same comparison arguments  yield  a contradiction with the aid of Hopf's boundary point lemma at $z^* \in \partial D$,  and hence case (IIc-1-2) does not occur. Eventually, case (IIc-1) does not occur.  We easily know that the same manner as in case (IIc-1) works also in case (IIc-3).  

Let us proceed to the remaining case (IIc-2). Here we need Proposition \ref{prop:ODE-PDE comparison general}. Distinguish two cases:
$$
\mbox{ (IIc-2-1) } \ \sigma_c > \sigma_s; \quad \mbox{ (IIc-2-2) } \ \sigma_c < \sigma_s.
$$
In case (IIc-2-2), we employ $v_{r_0}$. It follows from (3) of Proposition \ref{prop:ODE comparison with different conductivities} that 
$$
v_{r_0}(\rho) > w_+(\rho) \ \mbox{ for every } \rho \not=r_0,\ \mbox{ and hence } c_4(r_0) > 0.
$$
Because of \eqref{the boundary of D in z coordinates} there exists a point $z^0 \in \partial D$ with $z^0_N = r_0$ and moreover 
$$
v_{r_0} = w \ \mbox{ and }\ \nabla v_{r_0} = \nabla w = 0\ \mbox{ at the point } z^0 \in \partial D.
$$
Then  the same comparison arguments  yield  a contradiction with the aid of Hopf's boundary point lemma at $z^0 \in \partial D$. Thus, case (IIc-2-2) does not occur.

In case  (IIc-2-1), we employ $v_*$. It follows from (2) of Proposition \ref{prop:ODE comparison with different conductivities} that 
\begin{equation}
\label{inequality needed not complete}
w_+(\rho) < v_{*}(\rho)  \ \mbox{if }  z_N^*< \rho \le r_0.
\end{equation}
Remark that this inequality is not sufficient for the previous comparison arguments, because of \eqref{the boundary of D in z coordinates}.
For the sake of this reason, by integrating the ordinary differential equations which $w_+$ and $v_*$ satisfy, we have from \eqref{inequality needed not complete}
$$
\sigma_c v_*^\prime(r_0)= \sigma_c v_*^\prime(r_0) - \sigma_s w_+^\prime(r_0)= \int_{z_N^*}^{r_0}(v_*(\rho)-w_+(\rho)) d\rho > 0.
$$
Hence $v_*^\prime(r_0) > 0$.  By choosing a constant $\gamma > 0$ satisfying
$$
v_*(r_0) = \gamma \exp\left( - \frac {r_0}{\sqrt{\sigma_c}}\right),
$$
we introduce a function $v_{**} = v_{**}(\rho)$ for $\rho \ge 0$ given by
$$
v_{**}(\rho) = 
\begin{cases}
  \gamma \exp\left( - \frac {\rho}{\sqrt{\sigma_c}}\right) &\mbox{ if }\ r_0 \le \rho,
\\
v_*(\rho)\ &\mbox{ if }\  z_N^* \le \rho < r_0\,
\\
w_+(\rho)\ &\mbox{ if }\  0 \le \rho < z_N^*.
 \end{cases}
 $$
 Hence we have in particular
 \begin{equation}
\label{key inequality in Proposition 6.4}
 (\sigma_cv_{**}^\prime(\rho) - \sigma_s w_+^\prime(\rho))(v_{**}^\prime(\rho) - w_+^\prime(\rho)) > 0\ \mbox{ if } z_N*< \rho < z_N^\sharp.
\end{equation}
Indeed, for $z_N*< \rho \le r_0$, by integrating the ordinary differential equations which $w_+$ and $v_*$ satisfy, we have from \eqref{inequality needed not complete}
$$
\sigma_cv_{**}^\prime(\rho) - \sigma_s w_+^\prime(\rho) = \sigma_cv_{*}^\prime(\rho) - \sigma_s w_+^\prime(\rho) = \int_{z_N^*}^\rho(v_{*}(s)-w_+(s)) ds > 0.
$$
Then,  since $w_+^\prime(\rho) < 0$ and $\sigma_c > \sigma_s$, we have
$$
v_{**}^\prime(\rho) - w_+^\prime(\rho) =  \frac 1{\sigma_c}\left(\sigma_cv_{**}^\prime(\rho) - \sigma_c w_+^\prime(\rho)\right) > \frac 1{\sigma_c}\left(\sigma_cv_{**}^\prime(\rho) - \sigma_s w_+^\prime(\rho)\right) > 0.
$$
Therefore,  for $z_N*< \rho \le r_0$,  inequality \eqref{key inequality in Proposition 6.4} holds true.  For $r_0 < \rho < z_N^\sharp$, since $v_{**}^\prime(\rho) < 0$ and $w_+^\prime(\rho) > 0$, inequality \eqref{key inequality in Proposition 6.4} follows easily. Moreover, since $v_{**}^\prime(r_0-0) > 0 > v_{**}^\prime(r_0+0)$ and $v_{**}(r_0-0) = v_{**}(r_0+0)$, we see that
$$
-\left(\sigma_2 v_{**}^\prime\right)^\prime + v_{**} \ge 0\ \mbox{ in } (0,\infty),
$$
where we set
$$
\sigma_2 = \sigma_2(\rho) = 
\begin{cases}
\sigma_s\  &\mbox{ if }\ 0 \le \rho \le z_N^*,
\\
\sigma_c\  &\mbox{ if }\  \rho > z_N^*.
\end{cases}
$$
Then we can apply Proposition \ref{prop:ODE-PDE comparison general} to $w_1 = w,\ w_2= v_{**},\ \ell = z_N^*$ and $L = z_N^\sharp$ and conclude that
$$
w \le v_{**} \ \mbox{ in } \Omega,\ \mbox{ and hence } w < v_{**} \ \mbox{ in } D.
$$
Therefore,  this yields  a contradiction with the aid of Hopf's boundary point lemma at $z^* \in \partial D$, and  case  (IIc-2-1) does not occur.  The proof of Theorem 
\ref{th:elliptic overdetermined} is complete. \qed


\setcounter{equation}{0}
\setcounter{theorem}{0}

\def\theequation{A.\arabic{equation}}
\def\thetheorem{A.\arabic{theorem}}

\section{Appendices}
\label{section6}
We deal with three comparison principles and one maximum principle for elliptic equations with discontinuous conductivities.
We start with  a comparison principle for two solutions of ordinary differential equations with different conductivities
(see Lemma 3.5 in \cite{SBessatsu2017}).
\begin{proposition}
\label{prop:ODE comparison with different conductivities}
Let $\sigma_j\ (j=1,2)$ be two constants with $0 < \sigma_1 < \sigma_2$ and let $w_j = w_j(\rho)\ (j=1, 2)$ solve $-\sigma_j w_j^{\prime\prime} + w_j = 0$ in $\mathbb R$ for $j=1, 2$,  respectively.   Suppose that $w_1(r) = w_2(r)$ for some $r \in \mathbb R$. Then the following assertions hold:
\begin{enumerate}[\rm (1)]
\item Assume that $\sigma_1w_1^\prime(r) = \sigma_2w_2^\prime(r) >0$. Then we have
\begin{itemize}
\item[{\rm (i)}] If there exists $s \in (-\infty, r)$ such that  $w_1(s) = w_2(s)$ and $w_1(\rho) < w_2(\rho)$ for every $\rho \in (s, r)$, then $w_1^\prime(s) < 0$ and $w_2^\prime(s) < 0$.
 
 \item[{\rm (ii)}] If there exists $\ell \in (r, \infty)$ such that  $w_1(\ell) = w_2(\ell)$ and $w_1(\rho) > w_2(\rho)$ for every $\rho \in (r,\ell)$, then $w_1^\prime(\ell) < 0$ and $w_2^\prime(\ell) < 0$.
 \end{itemize}
\item Assume that $\sigma_1w_1^\prime(r) = \sigma_2w_2^\prime(r) <0$. Then we have

\begin{itemize}
\item[{\rm (i)}] If there exists $s \in (-\infty, r)$ such that  $w_1(s) = w_2(s)$ and $w_1(\rho) > w_2(\rho)$ for every $\rho \in (s, r)$, then $w_1^\prime(s) > 0$ and $w_2^\prime(s) > 0$.
 
 \item[{\rm (ii)}] If there exists $\ell \in (r, \infty)$ such that  $w_1(\ell) = w_2(\ell)$ and $w_1(\rho) < w_2(\rho)$ for every $\rho \in (r,\ell)$, then $w_1^\prime(\ell) > 0$ and $w_2^\prime(\ell) > 0$.
 \end{itemize}
 \item If $w_1^\prime(r) = w_2^\prime(r) =0$ and  $w_1(r) = w_2(r) > 0$, then  $w_1(\rho) > w_2(\rho)$ for every $\rho \not= r$.
\end{enumerate}
\end{proposition}

\noindent
{\it Proof.\ } Let us first consider (3).  Set $w_1(r) = w_2(r) = a >0$. Then it follows that for $j = 1, 2,$
$$
w_j(\rho) = \frac a2 \left\{ \exp\left(-\frac {\rho-r}{\sqrt{\sigma_j}}\right) +  \exp\left(\frac {\rho-r}{\sqrt{\sigma_j}}\right)\right\}\ \mbox{ for every } \rho \in \mathbb R.
$$
Since $0 < \sigma_1 < \sigma_2$, we have the conclusion.

Let us proceed to (1).  Note that
 \begin{equation}
 \label{equation of the difference for 3}
\sigma_1w_1^{\prime\prime}(\rho) - \sigma_2w_2^{\prime\prime}(\rho) = w_1(\rho)-w_2(\rho)\ \mbox{ for } \rho \in \mathbb R.
\end{equation}
Since $\sigma_1w_1^\prime(r) = \sigma_2w_2^\prime(r) > 0, \ w_1(r) = w_2(r)$ and $0 < \sigma_1 < \sigma_2$,  we observe that
$$
w_1^{\prime}(r) > w_2^{\prime}(r), 
$$
and hence there exists  a number $\delta > 0$ such that
 $$
 w_1(\rho) < w_2(\rho)\ \mbox{ for every } \rho \in (r-\delta, r)\ \mbox{ and  }\   w_1(\rho) >  w_2(\rho)\ \mbox{ for every } \rho \in(r, r + \delta).
 $$
Let us prove (i).  Since $\sigma_1w_1^\prime(r) = \sigma_2w_2^\prime(r),\  w_1(s) = w_2(s)$ and $w_1(\rho) < w_2(\rho)$ for every $\rho \in (s, r)$,  we notice that
$w_1^\prime(s)\le w_2^\prime(s)$.  Integrating \eqref{equation of the difference for 3} over the interval $[s, r]$ gives
$$
 - \sigma_1w_1^\prime(s) + \sigma_2w_2^\prime(s) = \int_s^r (w_1(\rho)-w_2(\rho))\ d\rho < 0.
$$
These yield that $w_1^\prime(s) < 0$ and $w_2^\prime(s) < 0$,  since $0 < \sigma_1 < \sigma_2$. (ii) is proved similarly.

It remains to consider (2).  Since $\sigma_1w_1^\prime(r) = \sigma_2w_2^\prime(r) < 0, \ w_1(r) = w_2(r)$ and $0 < \sigma_1 < \sigma_2$,  we observe that
$$
w_1^{\prime}(r) < w_2^{\prime}(r), 
$$
and hence there exists  a number $\delta > 0$ such that
 $$
 w_1(\rho) > w_2(\rho)\ \mbox{ for every } \rho \in (r-\delta, r)\ \mbox{ and  }\   w_1(\rho) <  w_2(\rho)\ \mbox{ for every } \rho \in(r, r + \delta).
 $$
 Thus the conclusion follows from the same argument as in (1). \qed

We have a proposition concerning the unique determination of  discontinuity of the conductivity for an ordinary differential equation with a nontrivial Cauchy data  (see Lemma 3.1 in \cite{SBessatsu2017} for the case dealing with bounded domains).
\begin{proposition}
\label{prop:ODE comparison with discontinuous conductivities}
Let $0 < r_1 \le r_2 < \infty$.  Define $\sigma_j = \sigma_j(\rho)\ (j=1, 2)$ for $\rho \ge 0$ by
$$
\sigma_j(\rho) = \begin{cases} \sigma_s\ &\mbox{ if } 0 \le \rho \le r_j, \\ 
  \sigma_c\  &\mbox{ if } r_j < \rho, 
  \end{cases}
$$
where $\sigma_c, \sigma_s$ are positive constants with $\sigma_c  \not= \sigma_s$.  
Let $w_j = w_j(\rho)\ (j=1, 2)$ be bounded solutions of $-(\sigma_j w_j^\prime)^\prime + w_j = 0$ in $[0, \infty)$ satisfying
\begin{eqnarray*}
&&w_1(0) =w_2(0),\ \ w_1^\prime(0) = w_2^\prime(0),
\\
&&\mbox{  and either }\  w_1(0) \not=0\ \mbox{ or }\ w_1^\prime(0) \not=0.
\end{eqnarray*}
Then  $r_1 = r_2$ and $w_1 \equiv w_2$ in $[0,\infty)$.
\end{proposition}

\noindent
{\it Proof.}  Since $w_j\ (j=1, 2)$ are bounded, we see that there exist two constants $c_j\ (j=1, 2)$ satisfying
$$
w_j(\rho) = c_j\exp\left(-\frac \rho{\sqrt{\sigma_c}}\right)\ \mbox{ for every } \rho \ge \rho_j \mbox{ and for } j =1, 2.
$$
Transmission conditions yield that $w_j\ (j=1, 2)$ are continuous on $[0,\infty)$  and 
$$
\sigma_s w_j^\prime(r_j-0) = \sigma_c w_j^\prime(r_j+0)\  \mbox{ for }  j =1, 2.
$$
Hence we have
\begin{eqnarray*}
&&\int_0^\infty w_1 w_2\ dx = \int_0^{r_1}( \sigma_s w_1^\prime)^\prime w_2\ dx + \int_{r_1}^\infty( \sigma_c w_1^\prime)^\prime w_2\ dx
\\
&&= - \sigma_s w_1^\prime(0) w_2(0) +\sigma_s w_1^\prime(r_1-0) w_2(r_1) - \sigma_c w_1^\prime(r_1+0) w_2(r_1) - \int_0^\infty \sigma_1 w_1^\prime w_2^\prime \ dx
\\
&&=  - \sigma_s w_1^\prime(0) w_2(0) - \int_0^\infty \sigma_1 w_1^\prime w_2^\prime \ dx.
\end{eqnarray*}
Thus we obtain
\begin{equation}
\label{1-2 with 1}
\int_0^\infty w_1 w_2\ dx = - \sigma_s w_1^\prime(0) w_2(0) - \int_0^\infty \sigma_1 w_1^\prime w_2^\prime\ dx.
\end{equation}
Changing the roles of $w_j\ (j=1, 2)$ yields that
\begin{equation}
\label{1-2 with 2}
\int_0^\infty w_1 w_2\ dx = - \sigma_s w_2^\prime(0) w_1(0) - \int_0^\infty \sigma_2 w_1^\prime w_2^\prime\ dx.
\end{equation}
In the same way we also have
\begin{eqnarray}
\int_0^\infty w_1^2\ dx &=& - \sigma_s w_1^\prime(0) w_1(0) - \int_0^\infty \sigma_1 (w_1^\prime)^2\ dx, \label{1-1 with 1}
\\
\int_0^\infty w_2^2\ dx &=& - \sigma_s w_2^\prime(0) w_2(0) - \int_0^\infty \sigma_2 (w_2^\prime)^2\ dx, \label{2-2 with 2}
\end{eqnarray}
Therefore by combing \eqref{1-2 with 1} and   \eqref{1-2 with 2} with the initial condition we obtain
\begin{equation}
\label{mixed product must vanish in the intersection} 
\int_{r_1}^{r_2} w_1^\prime w_2^\prime\ dx = 0,
\end{equation}
since $\sigma_s \not = \sigma_c$. Then  it follows from these equalities and the initial condition that
\begin{eqnarray*}
&&\int_0^\infty (w_1-w_2)^2\ dx
\\
&& = -\int_{r_1}^{r_2} (\sigma_c (w_1^\prime)^2 +\sigma_s(w_2^\prime)^2) dx - \int_0^{r_1} \sigma_s(w_1^\prime-w_2^\prime)^2  dx- \int_{r_2}^\infty \sigma_c(w_1^\prime-w_2^\prime)^2 dx
\\
&&\le 0,
\end{eqnarray*}
which yields that $w_1 \equiv w_2$. Moreover, since $w_1$ is not constant because of the initial condition,  it follows that $r_1=r_2$. \qed

Let us next give a maximum  principle for an elliptic equation in unbounded domains in $\mathbb R^N$, whose proof can be modified in proving the next key proposition.

\begin{proposition}
\label{prop:maximum principle on unbounded domains}
Let $D \subset \mathbb R^N$ be an unbounded domain, and let  $\sigma=\sigma(x)\ (x\in D)$ be 
a general conductivity satisfying
$$
0 < \mu \le \sigma(x) \le M\ \mbox{ for every } x \in \mathbb R^N,
$$
where $\mu, M$ are positive constants.  Assume that $w \in H^1_{loc}(D)\cap L^\infty(D)\cap C^0(\overline{D})$ satisfies
$$
-\mbox{\rm div}(\sigma \nabla w) + \lambda w \ge 0\ \mbox{ in } D\ \mbox{ and }\ w \ge 0\ \mbox{ on } \partial D
$$
for some constant $\lambda > 0$. Then $w \ge 0$ in $D$,  and moreover, either $w > 0$ in $D$ or $w \equiv 0$ in $D$.
\end{proposition}

\begin{remark}
When $D$ is bounded, this proposition is well known and holds true  for every $\lambda \ge 0$. However, when $D$ is unbounded, this proposition is not true for $\lambda = 0$. Indeed, a  counterexample is given in {\rm \cite[p. 37]{ABR2001Springer}}, where $N \ge 3,\ D = \{ x \in \mathbb R^N : |x| > 1 \},\ \sigma(x)\equiv 1$ and  $w(x) = |x|^{2-N}-1$.
\end{remark}

\noindent
{\it Proof of Proposition \ref{prop:maximum principle on unbounded domains}.} Define $v = v(x)$ by
$$
v(x) = e^{-\delta|x|} w(x)\ \mbox{ for } x \in \overline{D},
$$
where  $\delta > 0$ is a constant which will be chosen later.
Then $v \in H^1_{loc}(D)\cap L^\infty(D)\cap C^0(\overline{D})$ and moreover 
\begin{equation}
\label{limit zero at infinity}
\lim_{|x| \to \infty} v(x) = 0,
\end{equation}
since $w \in L^\infty(D)$.
For every $ \varepsilon > 0$,  we consider a nonnegative function
$$
\varphi(x) = \max\{ -\varepsilon - v(x), 0 \}\  \mbox{ for } x \in \overline{D}.
$$
Since $v \in H^1_{loc}(D)\cap L^\infty(D)\cap C^0(\overline{D})$ and $v \ge 0\ \mbox{ on } \partial D$, it follows from \eqref{limit zero at infinity} that $\varphi$ is compactly supported in $D$ and $\varphi \in H^1_0(D)$, and hence $e^{-2\delta|\cdot|}\varphi(\cdot)  \in   H^1_0(D)$. Therefore we obtain
\begin{eqnarray}
 0 &\le& \int\limits_D\left\{\sigma(x) \nabla w(x) \!\cdot \!\nabla\!\left(\varphi(x) e^{-2\delta|x|}\right) + \lambda w(x) \varphi(x)e^{-2\delta|x|}\right\} dx \nonumber
\\
&=& \int\limits_{D\cap\{v<-\varepsilon\}}\!\!\!\!\!\!\! \sigma e^{-\delta|x|}\left\{\left(\delta{v}\frac x{|x|} +\nabla{v}\right)\!\cdot\!\left(\nabla\varphi-2\delta\varphi\frac x{|x|}\right) + \frac\lambda\sigma{ v}\varphi \right\} dx. \label{the integral including varepsilon and delta}
\end{eqnarray}
Notice  that 
$$
\varphi(x) = 
\begin{cases}
-\varepsilon - v(x)   &\mbox{if }\ v(x) < -\varepsilon, \\
0  &\mbox{if }\ v(x) \ge -\varepsilon,
\end{cases} 
\quad\mbox{ and }\quad
\nabla \varphi(x) = 
 \begin{cases}
-\nabla v(x)   &\mbox{if }\ v(x) < -\varepsilon, \\
0  &\mbox{if }\ v(x) \ge -\varepsilon.
\end{cases} 
$$
By setting
$$
I = \sigma^{-1} e^{\delta |x|} \times \mbox{ the integrand of the integral \eqref{the integral including varepsilon and delta}},
$$
we have
\begin{eqnarray*}
I &=& -|\nabla v|^2-\frac\lambda\sigma v^2+2\delta^2v^2+\delta v\frac x{|x|} \cdot\nabla v + \varepsilon\left(2\delta^2 v+2\delta\frac x{|x|} \cdot\nabla v-\frac\lambda\sigma v\right)
\\
&\le & -\left\{ 1-\delta\left(\frac 12+\varepsilon\right)\right\}|\nabla v|^2 - \left\{ \frac\lambda\sigma\left(1-\frac\varepsilon 2\right) -\left(2\delta^2+\frac \delta 2\right) \right\} v^2 + \varepsilon\left(\frac \lambda {2\sigma}+\delta\right).
\end{eqnarray*}
Here we have used Cauchy's inequality $2ab \le a^2+b^2$ and the fact that $v< 0$ in the integrand of \eqref{the integral including varepsilon and delta}.  Therefore, since $0 < \mu \le \sigma(x)\le M$, we can choose $\delta > 0 $ sufficiently small to obtain that if $0 < \varepsilon < 1$ then 
$$
I \le -\frac 14\left( |\nabla v|^2 + \frac\lambda M v^2\right) + \varepsilon\left(\frac \lambda {2\mu}+\delta\right)
$$
and hence
$$
\mu\!\!\!\!\!\!\!\!\! \int\limits_{D\cap\{v<-\varepsilon\}}\!\!\!\!\!\!\! e^{-\delta|x|}\left( |\nabla v|^2 + \frac \lambda M v^2\right) dx \le M\varepsilon\left(\frac {2\lambda}\mu+4\delta\right)\int\limits_{D} e^{-\delta |x|} dx.
$$
By choosing a sequence $\{\varepsilon_n\}$ with $\varepsilon_n\downarrow 0$ as $n \to \infty$ and letting $n \to \infty$, we conclude that
$$
\int\limits_{D\cap\{v<0\}}\!\!\!\!\! e^{-\delta|x|}\left( |\nabla v|^2 + \frac\lambda M v^2\right) dx = 0
$$
and hence $v \ge 0$ in $D$.  Therefore $w \ge 0$ in $D$. Once this is shown, the last part follows from the strong maximum principle (see \cite[Theorem 8.19, pp. 198--199]{GT1983}). \qed

Finally, we give a comparison principle  for  two solutions of differential inequalities with different conductivities on a half-space of $\mathbb R^N$ (see Lemma 3.3 in \cite{SBessatsu2017} for the case dealing with bounded domains).
\begin{proposition}
\label{prop:ODE-PDE comparison general} For two numbers  $L > \ell > 0$, set
$$
\Omega = \{ z \in \mathbb R^N : z_N > 0 \},\ E =\{ z \in  \mathbb R^N : z_N > \ell \}\ \mbox{  and } F =\{ z \in  \mathbb R^N : z_N > L \}.
$$
Let $D \subset \mathbb R^N$ be a domain with $C^2$ boundary $\partial D$ satisfying that
$
F \subset D \subset E.
$
Let $\sigma_j = \sigma_j(z) \ (j=1, 2)$ be given by
$$
\sigma_1 =
\begin{cases}
\sigma_c \quad&\mbox{in } D, \\
\sigma_s \quad&\mbox{in } \Omega \setminus D,
\end{cases}
\qquad\sigma_2 =
\begin{cases}
\sigma_c \quad&\mbox{in } E, \\
\sigma_s \quad&\mbox{in } \Omega \setminus E,
\end{cases}
$$
where $\sigma_c, \sigma_s$ are positive constants with $\sigma_c  \not= \sigma_s$.   Let $w_1  \in H^1_{loc}(\Omega) \cap L^\infty(\Omega),\ w_2  \in H^1_{loc}((0,\infty)) \cap L^\infty( (0,\infty))$ satisfy
\begin{eqnarray*}
&& -\mbox{\rm div}( \sigma_1 \nabla w_1) + w_1 =  0\  \mbox{ in }\ \Omega,\quad   -( \sigma_2 w^\prime_2)^\prime + w_2 \ge 0\  \mbox{ in } (0,\infty), 
\\
&& \mbox{ and }\  w_1(z) \equiv w_2(z_N)\ \mbox{ for }\ z \in \Omega \setminus E.
\end{eqnarray*}
 Then,  if  
$$
\left(\sigma_c w_2^\prime(z_N)- \sigma_s \frac {\partial w_1}{\partial z_N}(z)\right)\left( w_2^\prime(z_N)-  \frac {\partial w_1}{\partial z_N}(z)\right) \ge 0 \ \mbox{ for } z  \in E \setminus\overline{D}, 
$$
we have that  $w_1(z) \le w_2(z_N)\  \mbox{ for } z \in \Omega$.
\end{proposition}

\noindent
{\it Proof.} We modify the proof of Proposition \ref{prop:maximum principle on unbounded domains}.   First of all, we extend $w_2$ for $z =(z_1, \dots, z_{N}) \in \Omega$ by $w_2(z) = w_2(z_N)$.    Introduce a function $\psi = \psi(t)\ ( t \in \mathbb R)$ by
$$
\psi(t) = 
 \begin{cases}
e^{-\delta(t-L)}  &\mbox{if }\  t > L, \\
1 &\mbox{if }\ t \le L,
\end{cases} 
$$
where $\delta > 0$ is a constant which will be chosen later.  
Then we define $v=v(z)$ by
$$
 v(z) = e^{-\delta|\hat z|}\psi(z_N)(w_2(z)-w_1(z))\ \mbox{ for } z \in \overline{\Omega},
$$
where $\hat z = (z_1, \dots, z_{N-1}) \in \mathbb R^{N-1}$. Note that $v= 0$ in $\Omega\setminus E$.
  If $0 < \varepsilon < 1$,  we set 
$$
\varphi(z) = \max\{ -\varepsilon - v(z), 0 \}\ ( \ge 0) \mbox{ for } z \in \Omega.
$$
Since $\varphi$ is compactly supported in $\Omega$, we notice that the function $\varphi(z) e^{-2\delta|\hat z|}\psi^2(z_N)$ belongs to $H_0^1(\Omega)$. Therefore we observe that
$$
0 \le \int\limits_{\Omega}\! \left\{ (\sigma_2(z) \nabla w_2(z)-\sigma_1(z)\nabla w_1(z))\!\cdot\! \nabla\!\!\left(\varphi(z) e^{-2\delta|\hat z|}\psi^2(z_N)\right) + v(z)\varphi(z) e^{-\delta|\hat z|}\psi(z_N)\right\} dz.
$$
Then, since $\varphi = 0$ in $\Omega\setminus E$,  we have
\begin{eqnarray}
&&0 \le  \int\limits_{E\setminus D}\!  (\sigma_c \nabla w_2(z)-\sigma_s\nabla w_1(z))\!\cdot\! \nabla\!\!\left(\varphi(z) e^{-2\delta|\hat z|}\psi^2(z_N)\right) dz\qquad\qquad\qquad\qquad\qquad\label{partwise 1}
\\
 && + \int\limits_{D}  \sigma_c \nabla(w_2(z)-w_1(z))\!\cdot \!\nabla\!\!\left(\varphi(z) e^{-2\delta|\hat z|}\psi^2(z_N)\right) dz +\!  \int\limits_{E}\!v(z)\varphi(z) e^{-\delta|\hat z|}\psi(z_N) dz.\qquad\quad\label{partwise 2}
 \end{eqnarray}
 By observing that $\psi=1$ and $w_1$ depends only on $z_N$  in $E\setminus D$, we see that the integral of \eqref{partwise 1} equals
 $$
- \!\!\!\!\!\!\!\!\!\!\!\!\int\limits_{(E\setminus D)\cap\{ v < -\varepsilon \}}\!\!\!\!  \left(\sigma_c \frac {\partial w_2}{\partial z_N}(z)-\sigma_s \frac {\partial w_1}{\partial z_N}(z)\right) \left(\frac{\partial w_2}{\partial z_N}(z)-\frac{\partial w_1}{\partial z_N}(z)\right) e^{-3\delta|\hat z|} dz\ ( \le 0).
 $$
As for the first integral of \eqref{partwise 2}, since we observe that
$$
\left|\nabla\left(e^{\delta|\hat z|}(\psi(z_N))^{-1}\right)\right| \le 2\delta e^{\delta|\hat z|}(\psi(z_N))^{-1}\ \mbox{ and } \left|\nabla\left(e^{-2\delta|\hat z|}\psi^2(z_N)\right)\right| \le 4\delta e^{-2\delta|\hat z|}\psi^2(z_N),
$$
the first integral of \eqref{partwise 2} is bounded from above by
$$
\sigma_c \!\!\!\!\!\!\int\limits_{D\cap\{ v < -\varepsilon \}} \!\!\!\!\!\! e^{-\delta|\hat z|}\psi(z_N)\left\{-|\nabla v|^2 +4\delta|\nabla v|\varphi + 2\delta |v||\nabla\varphi| +  8\delta^2|v|\varphi\right\} dz.
$$
Moreover, since $\varphi = -\varepsilon-v$,  with  Cauchy's inequality in hand,  we see that the first integral of \eqref{partwise 2} is bounded from above by
$$
\sigma_c \!\!\!\!\!\!\int\limits_{D\cap\{ v < -\varepsilon \}} \!\!\!\!\!\! e^{-\delta|\hat z|}\psi(z_N)\left\{-(1-3\delta)|\nabla v|^2 +\delta(3+8\delta)|v|^2\right\} dz.
$$
On the other hand,  since $0 < \varepsilon < 1$, the second integral of  \eqref{partwise 2} is bounded from above by
$$
\int\limits_{E\cap\{ v < -\varepsilon \}}\!e^{-\delta|\hat z|}\psi(z_N)(- \frac 12 |v|^2 + \frac 12\varepsilon) dz.
$$
Therefore,  in view of \eqref{partwise 1} and \eqref{partwise 2},  since $D \subset E$, we choose $\delta > 0 $ sufficiently small to conclude that if $0 < \varepsilon < 1$  then
$$
\sigma_c \!\!\!\!\!\!\int\limits_{D\cap\{ v < -\varepsilon \}} \!\!\!\!\!\! e^{-\delta|\hat z|}\psi(z_N)|\nabla v|^2 dz + \int\limits_{E\cap\{ v < -\varepsilon \}}\!e^{-\delta|\hat z|}\psi(z_N) |v|^2 dz \le 2\varepsilon \int\limits_{E}\!e^{-\delta|\hat z|}\psi(z_N)dz.
$$
By choosing a sequence $\{\varepsilon_n\}$ with $\varepsilon_n\downarrow 0$ as $n \to \infty$ and letting $n \to \infty$, we infer that
$$
\sigma_c \!\!\!\!\!\!\int\limits_{D\cap\{ v < 0 \}} \!\!\!\!\!\! e^{-\delta|\hat z|}\psi(z_N)|\nabla v|^2 dz + \int\limits_{E\cap\{ v < 0 \}}\!e^{-\delta|\hat z|}\psi(z_N) |v|^2 dz = 0
$$
and hence $v \ge 0$ in $E$,  which completes the proof. \qed

\end{document}